\documentclass[12pt]{JHEP3}

\usepackage{amsmath, epsfig}
\newtheorem{remark}{Remark}
\newtheorem{example}{Example}

\def\cee{{\relax\hbox{$\inbar\kern-.3em{\rm C}$}}}

\newcommand{\be}{\begin{equation}}
\newcommand{\ee}{\end{equation}}
\newcommand{\bea}{\begin{eqnarray}}
\newcommand{\eea}{\end{eqnarray}}
\newcommand{\bml}{\begin{mathletters}}
\newcommand{\eml}{\end{mathletters}}
\newtheorem{theorem}{Theorem}

\newtheorem{proposition}{Proposition}

\def\darr#1{\raise1.5ex\hbox{$\leftrightarrow$}\mkern-16.5mu #1}
\def\roughly#1{\raise.3ex\hbox{$#1$\kern-.75em\lower1ex\hbox{$\sim$}}}

\def\IB{\relax\hbox{$\inbar\kern-.3em{\rm B}$}}
\def\IC{\relax\hbox{$\inbar\kern-.3em{\rm C}$}}
\def\ID{\relax\hbox{$\inbar\kern-.3em{\rm D}$}}
\def\IE{\relax\hbox{$\inbar\kern-.3em{\rm E}$}}
\def\IF{\relax\hbox{$\inbar\kern-.3em{\rm F}$}}
\def\IG{\relax\hbox{$\inbar\kern-.3em{\rm G}$}}
\def\IGa{\relax\hbox{${\rm I}\kern-.18em\Gamma$}}
\def\IH{\relax{\rm I\kern-.18em H}}
\def\IK{\relax{\rm I\kern-.18em K}}
\def\IL{\relax{\rm I\kern-.18em L}}
\def\IP{\relax{\rm I\kern-.18em P}}
\def\IR{\relax{\rm I\kern-.18em R}}
\def\IZ{\relax\ifmmode\mathchoice{
\hbox{\cmss Z\kern-.4em Z}}{\hbox{\cmss Z\kern-.4em Z}}
{\lower.9pt\hbox{\cmsss Z\kern-.4em Z}} {\lower1.2pt\hbox{\cmsss
Z\kern-.4em Z}} \else{\cmss Z\kern-.4em Z}\fi}
\def\II{\relax{\rm I\kern-.18em I}}

\def\ee#1{{\rm erf}\left(#1\right)}

\def\Tr{{\rm Tr}}

\def\dim{{\rm dim}}

\def\inbar{\,\vrule height1.5ex width.4pt depth0pt}

\font\cmss=cmss10 \font\cmsss=cmss10 at 7pt

\def\lref{\begingroup\obeylines\lr@f}
\def\lr@f#1#2{\gdef#1{\ref#1{#2}}\endgroup\unskip}

\def\and{{a^\dagger_n}}

\def\math@note#1{\gdef\@eqnlabel{LAB: #1}}

\title{LOCALIZATION OF EQUIVARIANT COHOMOLOGY FOR COMPACT AND
NON-COMPACT GROUP ACTIONS}
\author{A. A. BYTSENKO \\
Departamento de F\'{\i}sica, Universidade Estadual de
Londrina\\
Caixa Postal 6001, Londrina-Paran\'a, Brazil\\
E-mail address: abyts@uel.br}
\author{M. LIBINE\\
Department of Mathematics and Statistics,
University of Massachusetts,\\
Lederle Research Tower, 710 North Pleasant Street, Amherst, MA 01003\\
E-mail address: matvei@math.umass.edu}
\author{F. L. WILLIAMS\\
Department of Mathematics and Statistics,
University of Massachusetts,\\
Lederle Research Tower, 710 North Pleasant Street, Amherst, MA 01003\\
E-mail address: williams@math.umass.edu}

\abstract{ We give a brief introduction to the Berline-Vergne
localization formula for the finite-dimensional setting and
indicate how the Duistermaat-Heckman formula is derived from it.
We consider applications of the localization formula when it is
specialized to a maximal dimensional co-adjoint orbit. In
particular, the case when the co-adjoint orbit is a quotient $G/T$
of a connected Lie group $G$ modulo a maximal torus $T$ is
analyzed in detail. We describe also a generalization of the
localization formula to non-compact group actions. }

\keywords{Equivariant cohomology; localization formula}
\begin{document}

\section{Introduction}

In 1982 J.J. Duistermaat and G. Heckman \cite{duistermaat} found a
formula which expressed
certain oscillatory integrals over a compact symplectic
manifold as a sum over
critical points of a corresponding phase function. In this sense these
integrals are localized, and their stationary-phase
approximation is exact
with no error terms occurring. The ideas and techniques of localization
extended to infinite-dimensional settings have proved to be
quite useful and
indeed central for many investigations in theoretical
physics - investigations
ranging from supersymmetric quantum mechanics, topological and
supersymmetric
field theories, to integrable models and low-dimensional gauge theories,
including two-dimensional Yang-Mills theory \cite{szabo}. Path integral
localization appears in the work of M. Semenov-Tjan-Schanskii
\cite{semenov}, which actually pre-dates \cite{duistermaat}.

E. Witten was the first to propose an extension of the
Duistermaat-Heckman (D-H) formula to an infinite-dimensional
manifold - namely to the loop space $LM$ of smooth maps from the
circle $S^1$ to a compact orientable manifold $M$. In this case a
purely formal application of the D-H formula to the partition
function of $N=1/2$ supersymmetric quantum mechanics yields a
correct formula for the index of a Dirac operator \cite{atiyah1}.
Further arguments in this direction were presented with
mathematical rigor by J.-M. Bismut in \cite{bismut1,bismut2}.

The various generalizations of D-H generally require formulations
in terms of equivariant cohomology. One has, for example, the
Berline-Vergne (B-V) localization formula
%[3-6]
\cite{berline1,berline2,berline3,berline4} which expresses the
integral of an equivariant cohomology class as a sum over zeros of
a vector field to which that class is related; also see
%[11-13]
\cite{woodhouse,blau,niemi1,szabo} for example, a
broader formulation of the localization formula.

This paper is organized as follows. In Sections 2-5 our remarks
are designed to provide the reader with a brief introduction to
the B-V localization formula, and to indicate how the D-H formula
is derived from it (see also \cite{bytsenko}). We limit our
discussion, in particular, to the finite-dimensional setting as
our idea is to convey the basic flavor of these formulas. This
introduction should prepare readers for quite more ambitions
discussions found in \cite{berline4,niemi1,szabo}, for example. In
Sections 6 we consider applications of the D-H localization
formula when it is specialized to a maximal dimensional co-adjoint
orbit. We pay attention to the case when the co-adjoint orbit is a
quotient $G/T$ of a connected Lie group $G$ (in particular the
unitary group $G=U(n)$) modulo a maximal torus $T$. Finally in
Section 7 we describe a generalization of the localization formula
to non-compact group actions.

The role of equivariant cohomology in physical theories will continue to
grow as it has grown in past years. In particular it will be an
indispensable tool for topological theories of gauge, strings, and
gravity.

\section{The equivariant cohomology space $H(M,X,s)$}

For an integer $j\geq 0$ let $\Lambda^jM$ denote the space of
smooth complex
differential forms of degree $j$ on a smooth manifold $M$.
$d:$ $\Lambda^jM\rightarrow \Lambda^{j+1}M$ will denote exterior
differentiation, and for a smooth vector field $X$ on $M$,
\begin{equation}
\theta(X): \Lambda^jM\rightarrow\Lambda^jM,\,\,\,\,\,
\iota(X): \Lambda^jM\rightarrow \Lambda^{j-1}M
\end{equation}
will denote Lie and interior differentiation by $X$, respectively:
\begin{eqnarray}
(\theta(X)\omega)(X_1,..., X_j) & = & X\omega(X_1,..., X_j)
\nonumber \\
& - & \sum_{\ell=1}^j\omega(X_1,..., X_{\ell-1},
[X,X_\ell], X_{\ell+1},..., X_j)\,, \\
(\iota(X)\omega)(X_1,..., X_{j-1}) & = & \omega(X,X_1,... X_{j-1})
\label{i0}
\end{eqnarray}
for $\omega\in \Lambda^jM$ and for
$X_1,..., X_j \in VM$ = the space of
smooth vector fields on $M$. One has the familiar rules
\begin{eqnarray}
\theta(X) & = & d\iota(X)+\iota(X)d\,,
\label{d1} \\
d\theta(X) & = & \theta(X)d, \,\,\, \theta(X)\iota(X)
=\iota(X)\theta(X)
\,,
\label{d2} \\
\iota(X)\circ \iota (X) & = & 0\,, \,\,\,\,\,d\circ d =0
\label{d3}
\mbox{.}
\end{eqnarray}
For a complex number $s$ let
\begin{equation}
d_{X,s}\stackrel{def}= d+s^{-1}\iota(X)\,\,\,{\rm on}\,\,\Lambda M
=\bigoplus_{j\geq0}\Lambda^jM
\mbox{.}
\label{d4}
\end{equation}
Then by (\ref{d1})--(\ref{d3}), $d_{X,s}\theta(X)=\theta(X)d_{X,s}$
and $d_{X,s}^2=s^{-1}\theta(X)$. Hence the subspace
\begin{equation}
\Lambda_XM=\{\omega \in \Lambda M | \theta(X)\omega=0\}
\mbox{,}
\label{s1}
\end{equation}
of $\Lambda M$ is $d_{X,s}-$ invariant and $d_{X,s}^2=0$
on $\Lambda_XM$. It follows that we can define the cohomology space
\begin{equation}
H(M,X,s)=Z(M,X,s)/B(M,X,s)
\label{H1}
\end{equation}
for $Z(M,X,s)=$ kernel of $d_{X,s}$ on $\Lambda_XM,\,\, \,
B(M,X,s)=d_{X,s}\Lambda_X M$. The space $H(M,X,s)$ appears
to depend on the
parameter $s$. However it is not difficult to show that for
$s\neq 0$ there is
an isomorphism of $H(M,X,s)$ onto $H(M,X,1)$.
For $X=0$, $H(M,0,s)$ is the ordinary de Rham cohomology of $M$.

We shall be interested in the case when $M$ has a smooth Riemannian
structure
$< , >$, and when $M$ is oriented and even-dimensional. Thus let
$\omega\in \Lambda^{2n}M-\{0\}$, ${\rm dim}\,M=2n$,
define the orientation of $M$. In this case we assume moreover
that $X$ is a Killing vector field:
\begin{equation}
X<X_1, X_2>\, = \,<[X, X_1], X_2>+<X_1, [X, X_2]>
\label{X1}
\end{equation}
for $X_1, X_2 \in VM$. If $p\in M$ is a zero of $X$ (i.e. $X_p=0$)
then there
is an induced linear map ${\frak L}_p(X)$ of the tangent space $T_p(M)$ of
$M$ at $p$ such that
\begin{equation}
{\frak L}_p(X)(Y_p)=[X, Y]_p\,\,\,\,\,{\rm for}\,\,\,Y\in VM
\mbox{.}
\label{L1}
\end{equation}
Because of (\ref{X1}) one has that ${\frak L}_p(X)$ is skew-symmetric;
i.e.
\begin{equation}
<{\frak L}_p(X)V_1, V_2>_p = -<V_1,{\frak L}_p(X)V_2>_p
\,\,\,\,\, {\rm for}\,\,\,\,\,
V_1,V_2 \in T_p(M)
\mbox{.}
\label{L2}
\end{equation}
Let
\begin{equation}
f_p(X): T_p(M)\oplus T_p(M)\rightarrow {\Bbb R}
\end{equation}
be the corresponding skew-symmetric bilinear form on $T_pM$:
\begin{equation}
f_p(X)(V_1,V_2) = <V_1,{\frak L}_p(X)V_2>_p\,\,\,
{\rm for}\,\,\,V_1,V_2\in T_pM
\mbox{.}
\label{f1}
\end{equation}
In order to apply some standard linear algebra to the real
inner product
space $(T_p(M),\, < , >_p)$, we suppose ${\frak L}_p(X)$ is a
non-singular linear
operator on $T_p(M): {\rm det}{\frak L}_p(X)\neq 0$; equivalently,
this means that the bilinear form $f_p(X)$ is non-degenerate.
Then one can find an ordered
orthonormal basis $e=e^{(p)}=\{e_j=e_j^{(p)}\}_{j=1}^{2n}$ of
$T_p(M)$ such that
\begin{eqnarray}
{\frak L}_p(X)e_{2j-1} & = & \lambda_je_{2j}\,,
\nonumber \\
{\frak L}_p(X)e_{2j} & = & -\lambda_je_{2j-1},\,\,\,\,\,{\rm for}
\,\,\,\,\, 1\leq j\leq n
\mbox{,}
\label{L3}
\end{eqnarray}
where each $\lambda_j\in {\Bbb R}-\{0\}$. In other words, relative to $e$
the matrix of ${\frak L}_p(X)$ has the form
\begin{equation}
{\frak L}_p(X) = \left[ \begin{array}{ccccccc}
0 & -\lambda_1 &  &  &  &  &    \\
\lambda_1 & 0  &  &  &  &  &    \\
            &  &  . &  &  &  &    \\
            &  &  & . &  &  &    \\
           &  &  &  & . &  &    \\
  &  &  &  &  &  0 & -\lambda_n \\
  &  &  &  &  &  \lambda_n & 0  \\
\end{array} \right]
\mbox{.}
\label{L4}
\end{equation}
Moreover, interchanging $e_1, e_2$ if necessary,
we can assume that $e$
is positively oriented: $\omega_p(e_1,..., e_{2n})>0$.
Finally, consider
the Pfaffian ${\rm Pf}_e({\frak L}_p(X))$ of ${\frak L}_p(X)$
relative to $e$:
\begin{equation}
{\rm Pf}_e({\frak L}_p(X))=\frac{1}{n!}\left[f_p(X)\wedge...\wedge
f_p(X)\right](e_1,...,e_{2n})
\mbox{.}
\label{Pf1}
\end{equation}
${\rm Pf}_e({\frak L}_p(X))$ satisfies
\begin{eqnarray}
{\rm Pf}_e({\frak L}_p(X))^2
& = & {\rm det}{\frak L}_p(X)\,,
\label{PF1}\\
{\rm Pf}_e({\frak L}_p(X))
& = &
(-1)^n\lambda_1\cdot\cdot\cdot \lambda_n
\mbox{.}
\label{Pf2}
\end{eqnarray}
If $e'=\{e_{j}'\}_{j=1}^{2n}$ is another ordered, positively oriented
orthogonal basis of $T_p(M)$ then
\begin{equation}
{\rm Pf}_{e'}({\frak L}_p(X))={\rm Pf}_e({\frak L}_p(X))
\mbox{.}
\label{Pf3}
\end{equation}
Equation (\ref{Pf3}) means that we can define a square-root of
${\frak L}_p(X)$ by setting
\begin{equation}
\left[{\rm det}{\frak L}_p(X)\right]^{1/2}=
(-1)^n{\rm Pf}_e({\frak L}_p(X))
\mbox{.}
\label{det1}
\end{equation}
That is, the square-root is independent of the choice $e$
of an ordered,
positively oriented orthogonal basis of $T_p(M)$.
By (\ref{Pf2}) we have
\begin{equation}
[{\rm det}{\frak L}_p(X)]^{1/2}=
\lambda_1\cdot\cdot\cdot \lambda_n
\mbox{.}
\label{det2}
\end{equation}
The reader is reminded that the hypotheses $X_p=0$ and
${\rm det}({\frak L}_p(X))\neq 0$ were imposed, with $X$
a Killing vector field.

\section{The localization formula}

As before we are given an oriented, $2n-$dimensional Riemannian manifold
$(M,\omega,< , >)$. Now assume that $G$ is a compact Lie group which acts
smoothly on $M$, say on the left, and that the metric $< , >$ is
$G-$invariant. Let ${\frak g}$ denote the Lie algebra of $G$. Given
$X\in {\frak g}$, there
is an induced vector field $X^{*}\in VM$ on $M$:
for $\phi\in C^{\infty}(M)$, $p\in M$
\begin{equation}
(X^{*} \varphi)(p)=\frac{d}{dt}\varphi(\exp(tX)\cdot p)|_{t=0}
\mbox{.}
\label{V1}
\end{equation}
Since $< , >$\, is $G-$invariant, one knows that $X^{*}$ is
a Killing vector
field. $X^{*}$ is said to be non-degenerate if, for every
zero $p\in M$ of
$X^{*}$, the induced linear map
${\frak L}_p(X^{*}): T_p(M)\rightarrow T_p(M)$ is
non-singular. Since $X^{*}$ is a Killing vector field,
${\frak L}_p(X^{*})$
is skew-symmetric with respect to the inner product structure
$< , >_p$ on
$T_p(M)$, as we have noted, and the non-singularity of ${\frak L}_p(X^{*})$
means that we can construct the square-root
\begin{equation}
\left[{\rm det}{\frak L}_p(X^{*})\right]^{1/2}=(-1)^n{\rm Pf}_e
({\frak L}_p(X^{*}))=\lambda_1\cdot\cdot\cdot \lambda_n
\mbox{,}
\label{det3}
\end{equation}
as in (\ref{det1}) and (\ref{det2}).

For a form $\tau\in \Lambda M =\bigoplus \Lambda^jM$ we write
$\tau_j\in \Lambda^jM$ for its homogeneous $j-$th component,
\begin{equation}
\tau = (\tau_0,..., \tau_{2n})=\sum_{j=0}^{2n}\tau_j
\mbox{,}
\label{t1}
\end{equation}
and we write $\left[{\tau}\right]$ for the cohomology class
of $\tau$ in case
$\tau\in Z(M,X,s)$ for $X\in VM,\, s\in {\Bbb C}$; i.e. $d_{X,s}\tau=0$
for $d_{X,s}$ in (\ref{d4}). When $M$ is compact, in particular,
one can integrate
any $2n-$form (as $M$ is orientable). Thus we can define
\begin{equation}
\int_M \tau = \int_M \tau_{2n}
\mbox{,}
\label{t2}
\end{equation}
and in fact we can define
\begin{equation}
\int_M [\tau] =\int_M \tau = \int_M \tau_{2n}
\mbox{.}
\label{t3}
\end{equation}
The integral $\int_M[\tau]$\, really does depend only on the
class $[\tau]$\, of $\tau$. Therefore the following result holds:

\begin{proposition}
If $\tau'\in B(M,X,s)$ then by a quick computation
using Stokes' theorem one sees that
\begin{equation}
\int_M\tau'= 0
\mbox{.}
\label{t4}
\end{equation}
Similarly if $p\in M$ with $X_p=0$ then $\tau_0'(p)=0$ for
$\tau'\in B(M,X,s)$.
\end{proposition}
\noindent \noindent {\bf Proof}. Indeed, if we write
$\tau'=d_{X,s}\beta$ for $\beta\in \Lambda_XM$ then one has
\begin{eqnarray}
\tau' && =  (s^{-1}\iota(X)\beta_1,\, d\beta_0+s^{-1}\iota(X)\beta_2,\,
d\beta_1+s^{-1}\iota(X)\beta_3,\,
d\beta_2+s^{-1}\iota(X)\beta_4,
\nonumber \\
&& ...,\, d\beta_{2n-2}+s^{-1}\iota(X)\beta_{2n},\, d\beta_{2n-1})
\nonumber \\
&&
= d\beta_0+s^{-1}\iota(X)\beta_0+d\beta_1+s^{-1}\iota(X)\beta_1+
d\beta_2+s^{-1}\iota(X)\beta_2 +
\nonumber \\
&&
...\, + d\beta_{2n}+s^{-1}\iota(X)\beta_{2n}
\mbox{.}
\label{t5}
\end{eqnarray}
Thus $\tau_0'(p)=[s^{-1}\iota(X)\beta_1]_{X=X_p}=
s^{-1}\beta_{1p}(X_p)=0$, and
\begin{equation}
\int_M\tau'=\int_Md\beta_{2n-1}=0
\mbox{,}
\label{t6}
\end{equation}
which proves (\ref{t4}).
$\Box$
\\

It follows that
the map $p^{*}: H(M,X,s)\rightarrow {\Bbb R}$ given by
\begin{equation}
p^{*}[\tau]= [\tau_0 \equiv s^{-1}\iota(X)\beta_1]_{X=X_p}
= \tau_0(p)\,\,\, {\rm for}\,\,\, X_p=0
\label{p1}
\end{equation}
is well-defined.
In
%[7-9],
\cite{berline1,berline2,berline3}, N. Berline and M. Vergne,
following some ideas of R. Bott in \cite{bott}, established the
following localization theorem, where the choice
$s^{-1}=-2\pi\sqrt{-1}$ is made.

\begin{theorem}
Assume as above that $M$ and $G$ are compact
and that the Riemannian metric $< , >$ on $M$ is $G-$invariant; i.e.
each $a\in G$ acts as an isometry of $M$. For $X\in {\frak g}$,
the Lie algebra
of $G$, assume that the induced vector field $X^{*}$ on $M$
(see (\ref{V1})) is
non-degenerate; thus the square-root in (\ref{det3}) is well-defined
(and is non-zero)
for $p\in M$ a zero of $X^{*}$ (i.e. $X_p^{*}=0$).
Then for any cohomology
class $[\tau]\in H(M,X^{*},s)$ one has
\begin{equation}
\int_M[\tau]=(-1)^{n/2}
\!\!\!\!\!
\sum_{\scriptstyle p\in M,
\atop\scriptstyle p=\,{\rm a\,\, zero\,\, of}\,\,X^{*}}
\frac{p^{*}[\tau]}
{[{\rm det}{\frak L}_p(X^{*})]^{1/2}}
\mbox{;}
\label{For1}
\end{equation}
see (\ref{t3}), (\ref{p1}).
\end{theorem}

For concrete applications of Theorem 1 we shall need to construct
concrete cohomology classes in $H(M,X^{*},s)$. The construction of
such classes requires that a bit more be assumed about $M$ and
$G$. Suppose for example that $M$ has a symplectic structure
$\sigma: \sigma \in \Lambda^2M$ is a closed two-form (i.e.
$d\sigma=0$) such that for every $p\in M$ the corresponding
skew-symmetric form
\begin{equation}
\sigma_p: T_p(M)\oplus T_p(M)\rightarrow {\Bbb R}
\end{equation}
is non-degenerate. In particular $M$ is oriented by the Liouville form
\begin{equation}
\omega_{\sigma}=
\frac{1}{n!}\sigma \wedge\cdot\cdot\cdot \wedge\sigma \in
\Lambda^{2n}M - \{0\}
\mbox{.}
\label{w1}
\end{equation}
Suppose also that there is a map
$J: {\frak g}\rightarrow C^{\infty}(M)$
which satisfies
\begin{equation}
\iota (X^{*})\sigma + dJ(X)=0,\,\,\,\,\,
\forall X\in {\frak g}
\mbox{,}
\label{i1}
\end{equation}
an equality of one-forms. The existence of such a map $J$
amounts to the
assumption that the action of $G$ on $M$ is Hamiltonian,
a point which we
shall return to later.

\begin{proposition}
For a given $J$ let us define for each
$X\in {\frak g}$ the form
$\tau^{X}\in \Lambda M$ by
\begin{equation}
\tau^{X}\stackrel{def}=
\left(J(X), 0, s\sigma, 0,..., 0
\right);
\label{t7}
\end{equation}
see (\ref{t1}). Then, we have that $\tau^{X}\in Z(M,X^{*},s)$.
\end{proposition}
\noindent \noindent {\bf Proof}. Since $J(X)$ is a function, then
$\iota (X^{*})J(X)=0$. Therefore by (\ref{d1})--(\ref{d3}) and
(\ref{i1}),
\begin{equation}
\theta(X^{*})J(X)=\iota (X^{*})dJ(X)=-\iota (X^{*})^2\sigma = 0
\label{th1}
\end{equation}
and
\begin{equation}
\theta(X^{*})\sigma = d\iota (X^{*})\sigma+\iota (X^{*})d\sigma
=d\iota (X^{*})\sigma\,\,\,\,\,
({\rm as}\,\,\, d\sigma=0) = -d^2J(X)=0
\mbox{.}
\label{th2}
\end{equation}
From the definition (\ref{t7}) it follows that
\begin{equation}
\theta(X^{*})\tau^X= \left(\theta(X^{*})J(X), 0,
\theta(X^{*})s\sigma,0,...,0\right)=0
\mbox{,}
\end{equation}
which by (\ref{s1}) means
that $\tau^X\in \Lambda_{X^{*}}M$. Also from the definitions
(\ref{d4}) and (\ref{i1}), we have
\begin{eqnarray}
d_{X^{*},s}\tau^X & = & (d+s^{-1}\iota (X^{*}))\tau^X
\nonumber \\
& = & dJ(X)+s^{-1}\iota (X^{*})J(X)+ds\sigma + s^{-1}\iota (X^{*})
s\sigma
\nonumber \\
& = & -\iota (X^{*})\sigma + \iota (X^{*})\sigma =0
\mbox{,}
\label{d5}
\end{eqnarray}
which verifies the claim, where again we have used that $\iota
(X^{*})J(X)= 0,\, d\sigma = 0$. $\Box$
\newline \newline Thus, for
a given $J$, we have for each $X\in {\frak g}$ a cohomology class
$[\tau^X] \in H(M,X^{*},s)$.

\section{The class{ [$e^{c\tau^X}$]}}

In the next Section the Duistermaat-Heckman formula will be derived by a
direct application of Theorem 1. The main point is the construction
of an appropriate cohomology class. Namely for the cocycle
$\tau^{X}\in Z(M,X^{*},s)$ in (\ref{t7}) we wish to
construct for $c\in {\Bbb C}$ a well-defined form
$e^{{c\tau}^X}$ which also is an element of
$Z(M,X^{*},s)$.

Thus again suppose that $J$ which satisfies (\ref{i1}) is given. For
$X\in {\frak g}$ let
\begin{equation}
\tau_0=J(X),\, \tau_1=0,\,\,\, \tau_2= s\sigma,\,\,\,
\tau_j=0\,\,\,\,\, {\rm for}\,\,\,\,\,
3\leq j \leq 2n
\mbox{,}
\label{t8}
\end{equation}
and let $\tau=\tau^X$. That is, by (\ref{t7}),
\begin{equation}
\tau =(\tau_0, \tau_1, \tau_2, ..., \tau_{2n})=(\tau_0, 0,
\tau_2, 0,0,...,0)
\label{t9}
\mbox{.}
\end{equation}
If $\omega_1, \omega_2$ are forms of degree $p,q$
respectively, then $\omega_1$ and $\omega_2$ commute if
either $p$ or $q$ is even, since
\begin{equation}
\omega_1\wedge \omega_2=(-1)^{pq}\omega_2\wedge \omega_1
\mbox{.}
\label{w2}
\end{equation}
In particular $\tau_0$ and $\tau_2$ commute. Now if $A$ and $B$
are commuting
matrices one has $e^{A+B}=e^A\cdot e^B$. Since $\tau_0$ and
$\tau_2$ commute
we should have, formally for any complex number $c$,
\begin{equation}
c\tau = c\tau_0+c\tau_2
\Rightarrow e^{c\tau}=e^{c\tau_0}\cdot e^{c\tau_2}
= e^{c\tau_0}(1+c\tau_2
+c^2\tau_2^2/2!+c^3\tau_2^3/3!+...)
\mbox{,}
\label{t10}
\end{equation}
with
\begin{equation}
\tau_2^j=\tau_2\wedge\cdot\cdot\cdot\wedge\tau_2
\,\,\,\,\,(j \,\,\,\,\,{\rm times})\,\,
\in \Lambda^{2j}M
\mbox{.}
\label{t11}
\end{equation}
Since $\Lambda^{2j}M=0$ for
$j>n$ we can take $\sum_{j=0}^{\infty}c^j\tau_2^j/j!$ to mean
$\sum_{j=0}^{n}c^j\tau_2^j/j!$. That is, thinking of $c\tau_2^j/j!$ as
$(0,0,...,c\tau_2^j/j!,0,...,0)$ and 1 as $(1,0,0,...,0)$ for
$1\in C^{\infty}(M)$, we are therefore lead to define $e^{c\tau}$ by
\begin{equation}
e^{c\tau} = \left(e^{c\tau_0},\,0\,,e^{c\tau_0}c\tau_2,\, 0\,,
e^{c\tau_0}
\frac{1}{2!}c^2\tau_2^2,\, 0,\, e^{c\tau_0}
\frac{1}{3!}c^3\tau_2^3,\, 0,
...,0, e^{c\tau_0}
\frac{1}{n!}c^n\tau_2^n\right) \in \Lambda M
\mbox{;}
\label{t12}
\end{equation}
which we can compare to expression (\ref{t1}).
Now $\iota (X^{*})e^{c\tau_0}=0$
(as $e^{c\tau_0}$ is a function),
and $de^{c\tau_0}=ce^{c\tau_0}d\tau_0$. That is, by
(\ref{d1})--(\ref{d3}),
\begin{eqnarray}
\theta(X^{*})e^{c\tau_0} & = &
c[\iota (X^{*})e^{c\tau_0}d\tau_0\,
+e^{c\tau_0}\iota (X^{*})d\tau_0]
\nonumber \\
& = & c\iota (X^{*})e^{c\tau_0}d\tau_0 =
ce^{c\tau_0}\iota (X^{*})d\tau_0
\mbox{,}
\label{t13}
\end{eqnarray}
since
\begin{eqnarray}
\tau_0 & = & J(X) \Rightarrow \,\,({\rm by}\,\,\, (\ref{d1})-(\ref{d3}),
\,\,(\ref{i1}))\,\,\,\,\,
\iota (X^*)d\tau_0 = -\iota (X^*)^2\sigma =0
\nonumber \\
& \Rightarrow & \theta(X^*)e^{c\tau_0}=0
\mbox{.}
\label{t14}
\end{eqnarray}
More generally,
\begin{eqnarray}
\theta(X^*)e^{c\tau_0}(c^j\tau_2^j)/j! & = &
(\theta(X^*)e^{c\tau_0})
(c^j\tau_2^j)/j! + e^{c\tau_0}(c^j/j!)\theta(X^*)\tau_2^j
\nonumber \\
& = & e^{c\tau_0}(c^j/j!)\theta(X^*)\tau_2^j\,\,\,
({\rm by}\,\,\,\, (\ref{t14})) = 0
\mbox{,}
\label{t15}
\end{eqnarray}
again by the fact that $\theta(X^*)$ is a
derivation and the fact that
$\theta(X^{*})\tau_2= s\theta(X^{*})\sigma$
with $\theta(X^{*})\sigma=0$
(as observed earlier).

\begin{proposition}
By (\ref{t12}) we see therefore that
\begin{equation}
\theta(X^{*})e^{c\tau}=0 \Rightarrow e^{c\tau}\in
\Lambda_{X^{*}}M
\mbox{,}
\label{t16}
\end{equation}
by (\ref{s1}). Therefore, we obtain
$d_{X^*,s}e^{c\tau}=0$.
\end{proposition}
\noindent \noindent By (\ref{t5}) and (\ref{t12})
\begin{eqnarray}
d_{X^*,s}e^{c\tau} & = & (0,\,d\beta_0+s^{-1}\iota (X^*)
\beta_2,\,0,\,d\beta_2+
s^{-1}\iota (X^*)\beta_4,\,0,
\nonumber \\
& ... & ,\,d\beta_{2n-2}+s^{-1}\iota(X^*)\beta_{2n},\,0 )
\label{d6}
\end{eqnarray}
for $\beta_{2j}=e^{c\tau_0}c^j\tau_2^j/j!$. Using that
\begin{equation}
d(\omega_1\wedge\omega_2) = d\omega_1\wedge\omega_2+
(-1)^{{\rm deg}\omega_1}\omega_1\wedge d\omega_2
\label{w3}
\end{equation}
for forms $\omega_1,\omega_2$
of homogeneous degree and that $e^{c\tau_0}, \tau_2$
are of even degree, we get
\begin{equation}
de^{c\tau_0}\tau_2^j=de^{c\tau_0}\wedge\tau_2^j +
e^{c\tau_0}\wedge d\tau_2^j
\mbox{,}
\label{d7}
\end{equation}
where
\begin{eqnarray}
d\tau_2^j & = & 0\,\,\,\, ({\rm by} \,\,\, (\ref{w3}))\,\,
\, {\rm since}\,\,\,
d\tau_2= sd\sigma=0
\nonumber \\
& \Rightarrow & d\beta_{2j}=(c^j/j!)e^{c\tau_0}dc\tau_0\wedge\tau_2^j
\nonumber \\
& = & -(c^{j+1}/j!)e^{c\tau_0}(\iota (X^*)\sigma)\wedge\tau_2^j
\mbox{,}
\label{t16}
\end{eqnarray}
by (\ref{i1}). Similarly
\begin{equation}
\iota(X^*)e^{c\tau_0}\tau_2^j  =  (\iota(X^*)e^{c\tau_0})\tau_2^j
+ e^{c\tau_0}\iota(X^*)\tau_2^j
= e^{c\tau_0}\iota(X^*)\tau_2^j
\mbox{,}
\label{i2}
\end{equation}
where
\begin{equation}
\iota(X^*)\tau_2^j = j\tau_2^{j-1}\wedge \iota(X^*)\tau_2
\end{equation}
since $\iota(X^*)$ also satisfies the derivative property (\ref{w3}),
and since $i(X^*)\tau_2$ and $\tau_2$ commute as
${\rm deg}\,\tau_2=2$. It follows
\begin{eqnarray}
j\tau_2^{j-1}\wedge \iota(X^*)\tau_2
& \Rightarrow & \iota(X^*)\beta_{2j}
e^{c\tau_0}(c^j/(j-1)!)\tau_2^{j-1}\wedge \iota(X^*)\tau_2
\nonumber \\
& = &
e^{c\tau_0}(c^j/(j-1)!)\tau_2^{j-1}\wedge \iota(X^*)s\sigma
\nonumber \\
& \Rightarrow &  s^{-1}\iota(X^*)\beta_{2j+2} =
e^{c\tau_0}(c^{j+1}/j!)\tau_2^j\wedge \iota(X^*)\sigma
\mbox{.}
\label{i3}
\end{eqnarray}
That is, by (\ref{t16}) and (\ref{i3}),
\begin{equation}
d\beta_{2j}+s^{-1}\iota(X^*)\beta_{2j+2}=0\,\, ({\rm again} \,\,\,\,
{\rm as} \,\,\,\, \iota (X^*)\tau_2
\,\,\,\, {\rm and} \,\,\,\,\tau_2\,\,\,\, {\rm commute})
\mbox{,}
\label{b1}
\end{equation}
which by (\ref{d6}) establishes our claim.
\\
\\
Hence the following result is holds:

\begin{theorem}
Suppose $J: {\frak g}\rightarrow
C^{\infty}(M)$
which satisfies (\ref{i1}) is given, where $\sigma$ is a symplectic
structure on
$M$. Recall that for $X\in {\frak g}$, equation (\ref{t7}) defines
a cocycle
$\tau^X\in Z(M,X^*, s)$. Similarly for $c\in {\Bbb C}$,
define $e^{c\tau^X}$ by (\ref{t12}):
\begin{equation}
e^{c\tau^X} =  \left(e^{cJ(X)}, 0, e^{cJ(X)}cs\sigma,
0, e^{cJ(X)}\frac{c^2}{2!}
(s\sigma)^2, 0,
... , 0, e^{cJ(X)}\frac{c^n}{n!}
(s\sigma)^n\right)\in \Lambda M
\mbox{,}
\label{T1}
\end{equation}
for ${\rm dim}\,M=2n$. Then also
$e^{c\tau^X}\in Z(M,X^*, s)$,
and thus we have the cohomology class $\left[e^{c\tau^X}\right] \in
H(M, X^*, s)$; see (\ref{d4}), (\ref{H1}), (\ref{V1}).
\end{theorem}

\section{The Duistermaat-Heckman formula}

Theorem 2 contains the basic assumption that a function $J: {\rm
{\bf g}}\rightarrow C^{\infty}(M)$ exists which satisfies
condition (\ref{i1}). As pointed out earlier this assumption
amounts to the assumption that the action of $G$ on $M$ is
Hamiltonian -- a point which we will now explain.

Given the symplectic structure $\sigma$ on $M$ there is a duality
$Y\leftrightarrow \beta_Y$ between smooth vector fields
$Y\in VM$ and smooth
one-forms $\beta_Y \in \Lambda^1M$ on $M$:
\begin{equation}
\beta_Y(X)=\sigma(Y,X)\,\,\,\,\,\,\,{\rm for\,\,\,every}\,\,\, X\in VM
\mbox{.}
\label{b2}
\end{equation}
$Y\in VM$ is called a Hamiltonian vector field if $\beta_Y$ is exact:
$\beta_Y = d\zeta$ for some $\zeta\in C^{\infty}(M)$. Let $HVM$
denote the
space of Hamiltonian vector fields on M. Actually $HVM$ is a Lie algebra.
For example, given any $\zeta\in C^{\infty}(M)$, the smooth
one-form $d\zeta$
corresponds (by the aforementioned duality) to a smooth vector
field $Y_\zeta$
on $M$. Thus $Y_\zeta\in HVM$ and by (\ref{i0}) and (\ref{b2})
we have for every $X\in VM$,
\begin{equation}
(\iota(Y_\zeta)\sigma)(X)=\sigma(Y_\zeta,X)=d\zeta(X)\, \Rightarrow
d\zeta=\iota (Y_\zeta)\sigma
\mbox{.}
\label{z1}
\end{equation}
The equation
\begin{equation}
[\zeta_1, \zeta_2]= Y_{\zeta_1}\zeta_2\,\,\,\,\,{\rm for}\,\,\,
\zeta_1,\zeta_2\in C^{\infty}(M)
\label{z2}
\end{equation}
defines the Poisson bracket $[\,,\,]$ on $C^{\infty}(M)$ which converts
$C^{\infty}(M)$ into a Lie algebra such that the map\,\,
$\wp : \zeta\rightarrow Y_\zeta: C^{\infty}(M)\rightarrow HVM$ is a Lie
algebra homomorphism; i.e.
\begin{equation}
[Y_{\zeta_1},Y_{\zeta_2}]=Y_{[\zeta_1,\zeta_2]}
\mbox{.}
\label{z3}
\end{equation}
The (left) action of $G$ on $M$ is called symplectic if
$X^*\in HVM,\,\,\, \forall X\in {\frak g}$; see (\ref{V1}).
Now the map $X\rightarrow X^*:
\,\,{\frak g}\rightarrow VM$ is not a Lie algebra homomorphism
since
\begin{equation}
[X_1,X_2]^*=-[X_1^*, X_2^*]\,\,\,\,\,\,\, {\rm for}
\,\,\,\,\, X_1, X_2 \in {\frak g}
\mbox{.}
\label{z4}
\end{equation}
If we define
\begin{equation}
\eta: {\frak g}\rightarrow VM\,\,\,\,\,\,\, {\rm by}
\,\,\,\,\, \eta(X)=(-X^*)=(-X)^*
\mbox{,}
\label{z5}
\end{equation}
then we do obtain a homomorphism:
\begin{equation}
\eta([X_1, X_2])=-[X_1, X_2]^* = [X_1^*, X_2^*]=
[-\eta(X_1), -\eta(X_2)] = [\eta(X_1), \eta(X_2)]
\mbox{.}
\label{z6}
\end{equation}
In other words if the action of $G$ is symplectic then
$\eta: {\frak g}\rightarrow HVM$ is a Lie
algebra homomorphism. The (left) action of $G$ on $M$ is called
Hamiltonian if it is symplectic and if the Lie algebra homomorphism
$\eta: {\frak g}\rightarrow HVM$ has a lift to
$C^{\infty}(M)$ -- i.e. if
there exists a Lie algebra homomorphism
$J: {\frak g}\rightarrow C^{\infty}(M)$
such that the diagram
\vspace{1.0cm}
\medskip
\par \noindent
$$
\begin{picture}(120,80)
\put(60,25){${\frak g}$} \put(0,80){${\rm C}^{\infty}({\rm M})$}
\put(100,80){${\rm HVM}$} \put(60,90){$\wp$}
\put(40,83){\vector(1,0){52}} \put(15,45){{\rm J}}
\put(50,40){\vector(-1,1){30}} \put(75,40){\vector(1,1){30}}
\put(100,45){$\eta$} \put(210,50) %{$(5.4)$}
{$ $}
\end{picture}
$$
is commutative: $\eta=\wp\circ J$, or
\begin{equation}
-X^* =Y_{J(X)}\,\,\,\,\,{\rm for\,\,\,every}\,\,\,X\in
{\frak g}
\label{z7}
\mbox{.}
\end{equation}
We note that such a $J$ will indeed satisfy condition (\ref{i1}).
Namely, by (\ref{z1}) and (\ref{z7}),
\begin{equation}
dJ(X)=\iota(Y_ {J(X)})\sigma=-\iota(X^*)\sigma\,\,\,\,\,\,\, {\rm for}
\,\,\,\,\, X\in {\bf {\frak g}}
\mbox{.}
\label{z8}
\end{equation}
The triple $(M,\sigma,J)$, for $J$ subject to (5.4), is called
a Hamiltonian
$G-$space \cite{kostant,woodhouse}. The basic example of a Hamiltonian
$G-$space is that
of an orbit ${\mathcal O}$ in the dual space ${\frak g}^*$ of
${\frak g}$ under the
co-adjoint action of $G$ on ${\frak g}^*$, which is induced
by the adjoint action of $G$ on ${\frak g}$, and where
$\sigma$ is chosen as the Kirillov
symplectic form on $M={\mathcal O}$, and $J$ is given by a
canonical construction.
Namely for a linear functional
$f$ on ${\frak g}$,\, $f\in {\frak g}^*$,
\begin{equation}
(a\cdot f)(X) = f(Ad(a^{-1})X)\,\,\,\,\,{\rm for}\,\,\,a\in G,\,\,
X\in {\frak g}
\mbox{.}
\label{a1}
\end{equation}
We shall recall how the (well-known) symplectic structure $\sigma$ on
${\mathcal O}$ is obtained (due to A.A. Kirillov)
and how the lifting $J$ is
canonically constructed. Thus we exhibit $({\mathcal O}, \sigma=
\sigma_{\mathcal O}, J=J_{\mathcal O})$ as a key example of a
Hamiltonian $G-$space. For this purpose
it is convenient to regard the orbit of $f$ as a homogeneous space:
${\mathcal O} \simeq G/G_f$ where $G_f$ is the stabilizer of $f$:
\begin{equation}
G_f = \{ a\in G | a\cdot f = f \}
\mbox{.}
\label{a2}
\end{equation}
$G_f$ is a closed subgroup of $G$ with Lie algebra ${\frak g}_f$
given by
\begin{equation}
{\frak g}_f = \{X \in {\frak g}\,\, |\, f([X, Y]) =0
\,\,\,\,\, \forall\, Y\in {\frak g} \}
\mbox{.}
\label{a3}
\end{equation}

Let $\tau^f$ be the corresponding Maurer-Cartan one-form
on $G$. That is,
$\tau^f \in \Lambda^1G$ is the unique left-invariant one-form on
$G$ subject to the condition
\begin{equation}
\tau^f(X)(1) = f(X)\,\,\,\,\,\,\,\,\,\,\, \forall \,X \in {\frak g}
\mbox{.}
\label{a4}
\end{equation}
Let $\pi: \,G\rightarrow G/G_f$ denote the quotient map.

\begin{theorem}  \label{sympl}
$G/G_f$ has a symplectic structure $\sigma$
which is uniquely given by $\pi^*\sigma = d\tau^f$.
\end{theorem}
\noindent

Here $\pi^*\omega_1$ denotes the pull-back of a form $\omega_1$.
The form $\sigma$ is also left-invariant; i.e. $\ell_a^*\sigma =
\sigma$ where $\ell_a: G/G_f \rightarrow  G/G_f$ denotes left
translation by $a\in G$. Given $X \in {\frak g}$ define $\psi_X:
G/G_f \rightarrow {\Bbb R}$ by
\begin{equation}
\psi_X(aG_f) = f(Ad(a^{-1})X) = (a\cdot f)(X)
\label{a5}
\end{equation}
for $a \in G$; $\psi_X$ is well-defined by (\ref{a2}).
One can show by computation that
\begin{equation}
d\psi_X = -\iota(X^*)\sigma
\mbox{.}
\label{a6}
\end{equation}
That is, by (\ref{b2}), $\beta_{-X^*} = d\psi_X \Rightarrow -X^*$
(or $X^*$) is Hamiltonian for each $X \in {\frak g}$; i.e. the
action of $G$ on $G/G_f$ is symplectic.
\\

\begin{theorem}
The action $G$ on $G/G_f$ is Hamiltonian.
\end{theorem}
\noindent
\noindent
{\bf Proof}. To see that this action is
Hamiltonian we must construct a lift $J: {\frak g}\rightarrow
C^{\infty}(G/G_f)$ of $\eta: X\rightarrow -X^*$. Namely define $J$
by
\begin{equation}
J(X) = \psi_X \,\,\,\,\,\, {\rm for}\,\,\,\,\psi_X\,\,\, {\rm in}\,\,\,
(\ref{a5})
\mbox{.}
\label{a7}
\end{equation}
Recall that $\wp: C^{\infty}(M) \rightarrow HVM$ is given by
$\wp(\psi) = Y_\psi$. That is, by (\ref{z1}) and (\ref{a6}),
$\wp(\psi_X) = - X^* = \eta (X)$, which shows that $J$ does
satisfy the commutative diagram (see above). The final step is to
show that $J$ is a homomorphism. Let $X_1, X_2 \in {\frak g}$, $a
\in G$. The Poisson bracket is given by (\ref{z2}):
\begin{eqnarray}
\left[J(X_1), J(X_2)\right](\pi(a)) & = &
\left(Y_{J(X_1)}J(X_2)\right)(\pi(a))
\nonumber \\
& = & \left(\wp(J(X_1)) J(X_2)\right)(\pi(a))
\nonumber \\
& = & \left(\eta(X_1)J(X_2)\right)(\pi(a))\,\,\,\,\,
({\rm again\,\,\, by}\,\,\, (5.4))
\nonumber \\
& = &
\left((-X_1^*)\psi_{X_2}\right)(\pi(a))
\,\,\,({\rm by}\,\,\,(\ref{a7}))
\nonumber \\
& = &
\frac{d}{dt}\psi_{X_2}\left((\exp (-tX_1))\cdot \pi(a)\right)|_{t=0}
\,\,\,\,\,({\rm by}\,\,\,\, (\ref{V1}))
\nonumber \\
& = &
\frac{d}{dt}\psi_{X_2}\left(\pi((\exp (-tX_1))\cdot a)\right)|_{t=0}
\nonumber \\
& = &
\frac{d}{dt}f\left(Ad(a^{-1}\exp(X_1))X_2\right)|_{t=0}
\,\,\,\,\,({\rm by}\,\,\,\, (\ref{a5}))
\nonumber \\
& = &
\frac{d}{dt}f\left(Ad(a^{-1})Ad(\exp(X_1))X_2\right)|_{t=0}
\nonumber \\
& = &
\frac{d}{dt}(a\cdot f)\left(Ad(\exp(X_1))X_2\right)|_{t=0}
\,\,\,\,\,({\rm by}\,\,\,\, (\ref{a5}))
\nonumber \\
& = &
(a\cdot f)\left([X_1, X_2]\right) = f\left(Ad(a^{-1})[X_1, X_2]\right)
\mbox{.}
\label{a8}
\end{eqnarray}
On the other hand
\begin{eqnarray}
J\left([X_1, X_2]\right)(\pi(a)) & = & \psi_{[X_1, X_2]}(\pi(a))
\,\,\,\,\, ({\rm by}\,\,\,(\ref{a7}))
\nonumber \\
& = &
f\left(Ad(a^{-1})[X_1, X_2]\right)\,\,\,\,\, ({\rm by}\,\,\,
(\ref{a5}))
\end{eqnarray}
which proves that $[J(X_1), J(X_2)] = J([X_1, X_2])$.
$\Box$
\\

We are now in position to state the Duistermaat-Heckman
formula -- in a form directly derivable from Theorem 1.

\begin{theorem}   \label{DHformula}
Suppose as above that $(M,\sigma,J)$ is a
Hamiltonian $G-$space where $G$ and $M$ are compact. Orient $M$ by the
Liouville form $\omega_\sigma$ in (\ref{w1}). Then for
$c\in {\Bbb C}$ and for
$X\in {\frak g}$ with $X^*$ non-degenerate, we have
\begin{equation}
\int_Me^{cJ(X)}\omega_\sigma = \left(\frac{2\pi}{c}\right)^n
\sum_{\scriptstyle p\in M,
\atop\scriptstyle p=\,{\rm a\,\, zero\,\, of}\,\,X^{*}}
\frac{e^{cJ(X)(p)}}{\left[{\rm det}\,{\frak L}_p(X^*)
\right]^{\frac{1}{2}}}
\mbox{.}
\label{MT}
\end{equation}
\end{theorem}

Here, as in Theorem 1, some $G-$invariant Riemannian metric $<,>$ on
$M$ has been selected, and the square-root in (\ref{MT}) is that
in (\ref{det3}).
\\
\\
{\bf Proof}.
The proof of (\ref{MT}) is quite simple, given Theorem 1.
Namely, given the
lifting $J$ (where we have noted that (5.4) implies (\ref{i1})) let
$c_J(X)=\left[e^{c\tau^X}\right]$ be the cohomology class constructed in
Theorem 2, for $c\in {\Bbb C},\,\, X\in {\frak g}$.
By (\ref{p1}) and (\ref{T1})
\begin{equation}
p^*c_J(X)=e^{cJ(X)(p)}\,\,\,\,\,{\rm for}\,\,\,X_p^*=0
\mbox{,}
\label{z9}
\end{equation}
and by (\ref{t3}) and (\ref{T1})
\begin{equation}
\int_Mc_J(X) = (sc)^n\int_Me^{cJ(X)}
\frac{\sigma^n}{n!}
= (sc)^n\int_Me^{cJ(X)}
\omega_\sigma
\mbox{.}
\label{z10}
\end{equation}
On the other hand given that $X^*$ is non-degenerate, the localization
formula (\ref{For1}) gives
\begin{equation}
\int_Mc_J(X) = (-1)^{\frac{n}{2}}
\sum_{\scriptstyle p\in M,
\atop\scriptstyle p=\,{\rm a\,\, zero\,\, of}\,\,X^{*}}
\frac{e^{cJ(X)(p)}}{\left[{\rm det}\,{\frak L}_p(X^*)
\right]^{\frac{1}{2}}}
\mbox{,}
\label{z11}
\end{equation}
by (\ref{z9}). That is, by (\ref{z10}) and (\ref{z11}) we obtain exactly
formula (\ref{MT}), as desired.
$\Box$
\\

Note that for $X\in {\frak g},\, Y\in VM$, and
$p\in M$,
\begin{eqnarray}
dJ(X)_p(Y_p) & = & [dJ(X)(Y)](p) =
[(-\iota(X^*)\sigma)(Y)](p)
\,\,\,({\rm as}\,\,\, J\,\,\, {\rm satisfies}\,\,\, (\ref{i1}))
\nonumber \\
& = &
-\sigma(X^*, Y)(p)\,\, ({\rm by} \,\,\,(\ref{i0}))\,\,\,=
-\sigma_p(X_p^*, Y_p)
\mbox{.}
\label{z12}
\end{eqnarray}
Hence $dJ(X)_p=0$ if $X_p^* = 0$, and conversely $dJ(X)_p = 0$\,\,\,
$\Rightarrow\,\,\,X_p^* = 0$ since $\sigma_p$ is non-degenerate.
(\ref{MT}) can therefore be expressed as
\begin{equation}
\int_Me^{cJ(X)}\omega_\sigma = \left(\frac{2\pi}{c}\right)^n
\sum_{\scriptstyle p\in M,
\atop\scriptstyle p=\,{\rm a\,\, critical\,\,point\,\, of}\,\,J(X)}
\frac{e^{cJ(X)(p)}}{\left[{\rm det}\,{\frak L}_p(X^*)
\right]^{\frac{1}{2}}}
\mbox{,}
\label{z13}
\end{equation}
where the critical points of $J(X)$ are those where $dJ(X)$ vanishes.
\\
\\
Recall that the asymptotic behaviour of an oscillatory integral
$
{\mathcal I}(f, t)\! =\! \int_{M}\! e^{\sqrt{-1}tf(x)}dx,
$
$
M =
$
some space, for large $t$ is given by the stationary-phase
approximation - the dominant
terms of this approximation being governed by the critical
points of the phase
$f(x)$. If we choose $c=\sqrt{-1}t$, for $t\in {\Bbb R}$,
in (\ref{z13}), in
particular, we see that the D-H formula can be viewed as an
exactness result in a stationary-phase approximation of the integrals
$\int_Me^ {\sqrt{-1}tJ(X)}\omega_\sigma$, as our remarks of Section 1
indicated.
For extended and much broader discussions of material introduced
here, the two
references \cite{berline4,szabo} are especially recommended.
The reference
\cite{szabo} in particular serves as a vast source of information
for the
needs of physicists. Further reading of interest is found in the
references
%[16-29].
\cite{kostant,atiyah2,witten1,stone,picken,niemi2,witten2,niemi3,witten3,dykstra,
niemi4,funahashi,cordes,paradan,schwarz}.

% Floyd's section
\section{Harish-Chandra, Itzykson-Zuber integral formulas}

Some very practical and beautiful applications of the general
D-H localization formula (Theorem \ref{DHformula}) result when
it is specialized, for example, to maximal dimensional co-adjoint
orbits ${\mathcal O}$.
To be specific, of special interest is the case where
${\mathcal O}$ is a quotient $G/T$ of a compact, connected Lie group $G$
modulo a maximal torus $T$.
For physically important reasons one often concentrates on the case of
the unitary group $G=U(n)$.
Integration over matrix groups $G$, which amounts to integration over
$G/T$
for $T-$invariant functions, has well-known importance for diverse areas
as quantum gravity \cite{FK}, integrable systems
\cite{SchlittgenWettig},
quantum chromodynamics \cite{FSV}, etc.
The Itzykson-Zuber (I-Z) integration formula (integration over $U(n)$)
\cite{ItzyksonZuber}, for example, occurs crucially in matrix models
(the Ising model on a random surface) where one considers the coupling
of conformal matter with two-dimensional quantum gravity.
This formula also appears in work on higher-dimensional lattice
gauge theories.
Harish-Chandra-Itzykson-Zuber integration over the symplectic
group $G$ is involved in
the computation of the mean of products of characteristic polynomials of random
matrices in certain ensembles \cite{FSV}.
Specialization of the D-H formula also leads to the Kirillov integral formula
for irreducible representations of $G$, which has relevance for geometric
quantization theory.

We establish in this section therefore the useful reductions of
Theorem \ref{DHformula} in the case of $G/T$, and we carry out
calculations
in the special, but important, case $G=U(n)$.
As a new localization formula will be presented in the next section
(due to the second named author) for a {\em non-compact} group,
it is helpful
here to gain further understanding of localization in the compact case.

$G$ will denote a compact, connected Lie group with Lie
algebra $\mathfrak g$
on which $G$ acts via the adjoint representation $Ad$.
Fix a Cartan subalgebra $\mathfrak t$ of $\mathfrak g$
(i.e. $\mathfrak t$ is a maximal abelian subalgebra) and
an $Ad(G)$-invariant innner product $\langle \cdot,\cdot \rangle$ on
$\mathfrak g$.
For $\mathfrak z$ the center of $\mathfrak g$, one has decompositions
\begin{equation}
\mathfrak g= \mathfrak z \oplus \mathfrak g_1= \mathfrak t \oplus
\mathfrak m,
\qquad \mathfrak t = \mathfrak z \oplus \mathfrak t_1
\end{equation}
for $\mathfrak g_1= [\mathfrak g,\mathfrak g]$
the commutator of $\mathfrak g$,
$\mathfrak t_1= \mathfrak t \cap \mathfrak g_1$,
$\mathfrak m=[ \mathfrak t_1, \mathfrak g_1 ]$.
If $T=\exp \mathfrak t$, which is a maximal torus in $G$ with
Lie algebra
$\mathfrak t$, then $\mathfrak m$ is $Ad(T)-$invariant which means that
the restriction of $\langle \cdot,\cdot \rangle$ to $\mathfrak m$
induces a natural $G-$invariant Riemannian metric
$\langle \cdot,\cdot \rangle$ on $M{\stackrel{\rm def}{=}} G/T$,
where one can regard $\mathfrak m$ as the tangent space to $M$ at
its origin 0.
An element $x\in \mathfrak g$ is called a {\em regular} element if its
centralizer
$\mathfrak g_x {\stackrel{\rm def}{=}} \{ y \in \mathfrak g |
\: [y,x]=0\}$
is of minimal dimension among all other centralizers:
$\dim \mathfrak g_x \le \dim \mathfrak g_y$ for all $y \in \mathfrak g$.
It is known that one can choose a regular element $x_0\in \mathfrak t_1$
such that $\mathfrak t= \mathfrak g_{x_0}$, $T=
G_{f_{x_0}}$ for $G_{f_{x_0}}$
given in (\ref{a2}) where $f_{x_0} \in \mathfrak g^*$
(the dual space of $\mathfrak g$) is given by
$f_{x_0}(x)=\langle x,x_0 \rangle$ for $x \in \mathfrak g$.
Thus $M$ is a co-adjoint orbit ${\mathcal O}$ of maximal dimension,
since $\mathfrak g_{x_0}$ is of minimal dimension.

The key ingredient in the D-H formula (\ref{MT}) is the determinant
there and its square root.
In order to describe these in the present setting a bit more notation is
necessary.
Choose a system of positive roots $P$ contained in the roots
$\Delta(\mathfrak g^{\mathbb C},\mathfrak t^{\mathbb C})$ of
$(\mathfrak g^{\mathbb C}, \mathfrak t^{\mathbb C})$,
the complexifications of $(\mathfrak g,\mathfrak t)$.
$\mathfrak g_1^{\mathbb C}$ is semisimple,
$\mathfrak t_1^{\mathbb C}$ is a Cartan subalgebra of
$\mathfrak g_1^{\mathbb C}$, and $P$ is related to a system
of positive roots
$\Delta ^+$ in the root system
$\Delta=\Delta(\mathfrak g_1^{\mathbb C}, \mathfrak t_1^{\mathbb C})$ of
$(\mathfrak g_1^{\mathbb C}, \mathfrak t_1^{\mathbb C})$ by a bijection
$\alpha \leftrightarrow \widetilde{\alpha}: \Delta^+
\leftrightarrow P$ where
$\widetilde{\alpha} (Z+x) {\stackrel{\rm def}{=}} \alpha (x)$ for
$Z+x \in \mathfrak z^{\mathbb C} \oplus \mathfrak t_1 ^{\mathbb C}
= \mathfrak t^{\mathbb C}$.
The main point here is that one can choose an orthonormal basis
$\mathcal{B} = \{ e_\alpha, f_\alpha \}_{\alpha \in \Delta ^+}$
of $\mathfrak m$ such that
\begin{equation}
[ H , e_\alpha ]= -\sqrt{-1} \alpha (H)f_\alpha,  \qquad
[ H , f_\alpha ]= \sqrt{-1} \alpha (H)e_\alpha \qquad
\text{for } H\in \mathfrak t_1
\end{equation}
(by compact Lie group structure theory), so that if a particular
ordering
$\{\alpha\}_{j=1} ^{n=\frac{1}{2}\dim M}$ of $\Delta^+$ is picked then
the matrix of $ad_H: \mathfrak m \rightarrow \mathfrak m$ relative to
$\mathcal{B}$ assumes the form
\begin{equation}  \label{6.1}
ad_H =
\begin{pmatrix}
0             & \sqrt{-1}\alpha_1(H) &         &               &              \\
-\sqrt{-1}\alpha_1(H) & 0            &         &               &              \\
              &              & \ddots  &               &              \\
              &              &         & 0             & \sqrt{-1}\alpha_n(H) \\
              &              &         & -\sqrt{-1}\alpha_n(H) & 0

\end{pmatrix},
\end{equation}
as in equation (\ref{L1}), for $H \in \mathfrak t_1$.
Moreover an element $X \in \mathfrak t$ is regular
$\iff$ $\alpha(X) \ne 0$
for all $\alpha \in \Delta$; and for a regular element
$X \in \mathfrak t_1$,
and $p \in M$ with $X^*_p = 0$ (see (\ref{V1})), say
$p=\pi(a)$ for $a \in G$,
where $\pi : G \to M=G/T$ is the natural map,
$\mathfrak{L}_p (X^*)$ in (\ref{s1}) considered as a map from
$\mathfrak m$ to $\mathfrak m$ is calculated to be given by
$\mathfrak{L}_p (X^*)=-ad_{Ad (a^{-1})X}$ where
$X^*_p=0 \Longrightarrow Ad (a^{-1})X \in \mathfrak t_1$.
By (\ref{6.1}) therefore, and the discussion in Sections 2, 3
one sees that
$X^*$ is nondegenerate and one can make the following choice of
square root of
$\det\mathfrak{L}_p(X^*)$, again for $p=\pi(a)$, $a\in G$:
\begin{equation}  \label{6.2}
\left[\det\mathfrak{L}_p(X^*)\right]^{\frac{1}{2}}
=\prod_{\alpha \in \Delta^+} \alpha (\sqrt{-1} Ad (a^{-1})X)
\mbox{,}
\end{equation}
where we assume that $\mathcal{B}$ is positively oriented
with respect to the
form $(\omega_\sigma)_0$ where $\omega_\sigma$ is the Liouville form
\begin{equation}
\omega_\sigma = \frac{\stackrel{\rm n-times}
{\overbrace{\sigma \wedge \ldots \wedge \sigma}}}{n!}
\end{equation}
in (\ref{w1}), $\sigma$ being the symplectic form on the orbit
${\mathcal O}=M =G/T=G/G_{f_{x_0}}$ and (again) 0 being the origin
$T$ of $M$. The latter assumption is satisfied if $x_0$ satisfies
\begin{equation}  \label{6.3}
\sqrt{-1}\alpha(x_0)>0 \,, \qquad \forall \alpha \in \Delta^+.
\end{equation}

The final point to make here is that (by computation) $X^*_{p=\pi(a)}=0$
(for any regular element $X \in \mathfrak t$) if and only if
\begin{equation}
a \in N_G(\mathfrak t) \quad {\stackrel{\rm def}{=}} \quad
\{ a \in G |\: Ad(a) \mathfrak t = \mathfrak t\},
\end{equation}
the normalizer of $\mathfrak t$ in $G$.
Here $N_G(\mathfrak t)$ is also the normalizer
\begin{equation}
N_G(T) \quad {\stackrel{\rm def}{=}} \quad
\{a \in G \, | \, aTa^{-1} = T \}
\end{equation}
of $T$ in $G$. On the other hand, the Weyl group of $(G,T)$ is
$W {\stackrel{\rm def}{=}} N_G(T)/T$ and one sees therefore that the map
$W \to M$ given by $w \mapsto p=\pi(a)$ for a coset
$w=aT \in W$, $a\in N_G(T)$, is a well-defined bijection of
$W$ onto the set
$F^X\, {\stackrel{\rm def}{=}} \, \{ p\in M \,|\, X^*_p =0 \}$,
which is the set that one sums over in (\ref{MT}).
Recalling that $T=G_{f_{x_0}}$, one has that the Hamiltonian lift
$J: g \to {\mathcal C}^\infty(M)$, given in general by (\ref{a5}),
(\ref{a7}), is given in the present situation by
\begin{equation}
J(Y)(aT) \quad {\stackrel{\rm def}{=}} \quad f_{x_0} (Ad (a^{-1})Y)
\quad {\stackrel{\rm def}{=}} \quad \langle Ad (a^{-1})Y, x_0 \rangle
\,\,\,\,\,\,\,  \text{for } a\in G.
\end{equation}
Putting the pieces together we see that the D-H localization formula
(\ref{MT})
reduces to the following concrete formula for $M=G/T$,
$c \in {\mathbb C}$, $\dim M=2n$:
\begin{eqnarray}  \label{6.4}
\int_{G/T} e^{ c \langle Ad (a^{-1})X,x_0 \rangle } \cdot
\frac{\stackrel{\rm n-times}
{\overbrace{\sigma \wedge \ldots \wedge \sigma}}}
{(2\pi)^n n!}(aT)
& = & c^{-n}\!\!\!\!\!\! \sum _{w\in W,\, w=aT,\, a \in N_G(T)}
\frac{e^{ c \langle Ad (a^{-1})X,x_0 \rangle }}
{\prod_{\alpha \in \Delta^+}\alpha(\sqrt{-1} Ad (a^{-1})X)}
\nonumber \\
& = & c^{-n}\sum_{w\in W}\frac{e^{c(w\cdot \lambda_{x_0})(X)}}
{\prod_{\alpha \in \Delta^+}(w\cdot \alpha)(\sqrt{-1} X)}
\mbox{,}
\end{eqnarray}
where $\lambda_{x_0} \in \mathfrak t_1^*$ is given by
$\lambda_{x_0}(H){\stackrel{\rm def}{=}} \langle H, x_0 \rangle$ for
$H \in \mathfrak t_1$ (i.e. $\lambda_{x_0}=f_{x_0}|_{\mathfrak t_1}$)
and
\begin{equation}
(w\cdot \lambda)(H) \quad {\stackrel{\rm def}{=}}
\quad \lambda(Ad (a^{-1})H)
\qquad \text{for } w=aT, \quad a \in N_G(T).
\end{equation}
Here we also assume that (as in (\ref{6.3})) $x_0$
satisfies the positivity
condition $\sqrt{-1} \alpha(x_0)>0$ for all $\alpha \in \Delta^+$.

We have noted that $y \in \mathfrak t_1$ is regular
$\iff \alpha(y) \ne 0$ for all $\alpha \in \Delta^+$,
and that $\alpha \leftrightarrow \widetilde{\alpha}$
is a bijection of $\Delta^+$ and $P$.
Since for $Y=Z+y \in \mathfrak z \oplus \mathfrak t_1= \mathfrak t$
one has that $\mathfrak g_{Y}=\mathfrak g_{y}$,
it follows that $Y$ is regular $\Longleftrightarrow y$ is regular
$\Longleftrightarrow \beta (Y) \ne 0$ for all $\beta \in P$,
and the localization formula in (\ref{6.4}) {\em extends} directly
from regular elements $X, x_0$ in $\mathfrak t_1$ to regular elements
$X,x_0$ in $\mathfrak t= \mathfrak z \oplus \mathfrak t_1$: for
$\omega {\stackrel{\rm def}{=}}
[(2\pi)^n n!]^{-1}\stackrel{\rm n-times}
{\overbrace{\sigma \wedge \ldots \wedge \sigma}}$,
\begin{equation}  \label{6.5}
\int_{G/T} e^{c \langle Ad (a^{-1})X,x_0 \rangle }\omega(aT)
=c^{-n} \sum_{w\in W}\frac{e^{c(w\cdot \lambda_{x_0})(X)}}
{\prod_{\beta \in P}(w\cdot \beta)(\sqrt{-1} X)}
\end{equation}
for $x_0$ satisfying $\sqrt{-1}\beta(x_0)>0$ for all
$\beta \in P$, where
$\lambda_{x_0} \in \mathfrak t^*$ is given by
$\lambda_{x_0}(H) {\stackrel{\rm def}{=}} \\
\langle H,x_0 \rangle$ for
$H \in \mathfrak t$ and
$(w \cdot \lambda)(Y) {\stackrel{\rm def}{=}} \lambda( Ad (a^{-1})Y)$
for $Y \in \mathfrak t^{\mathbb C}$, $w=aT \in W$, and $a \in N_G(T)$.
The symplectic structure $\sigma$ on $G/T$ is given by
Theorem \ref{sympl}.
Thus in (\ref{6.5}) we have arrived at Harish-Chandra's integral formula
\cite{HarishChandra, berline3, berline4}, which in essence computes the
Fourier transform of the measure $\omega$.

For the unitary group $G=U(n)$ with Lie algebra
$\mathfrak g= \mathfrak u(n)=$
the space of skew Hermitian matrices of degree $n$, one has
the following data:
$\langle X,Y \rangle {\stackrel{\rm def}{=}} -\Tr \,XY$ for
$X,Y \in \mathfrak g$,
$\mathfrak t=$ the space of diagonal matrices with entries
$\sqrt{-1}\theta_1,\ldots, \sqrt{-1}\theta_n$ for the
$\theta_j \in \mathbb{R}$,
$T=$ the group of diagonal matrices with entries
$e^{\sqrt{-1}\theta_1}, \ldots, e^{\sqrt{-1}\theta_n}$,
$\mathfrak g_1= \mathfrak{su}(n)$,
$\mathfrak t_1=$ matrices in $\mathfrak t$ with zero trace,
$\mathfrak z=$ matrices in $\mathfrak t$ with all entries equal,
$\mathfrak g^{\mathbb C}= \mathfrak{gl}(n,{\mathbb C})$,
$\mathfrak g_1^{\mathbb C}= \mathfrak{sl}(n,{\mathbb C})$,
$\Delta(\mathfrak g^{\mathbb C}, \mathfrak t^{\mathbb C})
=\{\widetilde{\alpha}_{rs}\}_{r\neq s}$,
where $\widetilde{\alpha}_{rs}(H)=H_r-H_s$ for
${\rm diag}(H_1, ..., H_n)$
%$H=\begin{pmatrix}
%H_1 & & \\
% & \ddots & \\
% & & H_n
%\end{pmatrix} \: \in \mathfrak t^{\mathbb C}=$
the space of complex diagonal matrices,
$\Delta(\mathfrak g_1^{\mathbb C},\mathfrak t_1^{\mathbb C})=$
the set of restrictions of elements of
$\Delta(\mathfrak g^{\mathbb C},\mathfrak t^{\mathbb C})$
to the trace zero matrices $\mathfrak t_1^{\mathbb C}$ in
$\mathfrak t^{\mathbb C}$.
\begin{equation}
P=\{\widetilde{\alpha}_{rs} | \: 1 \le r < s \le n \},  \qquad
\Delta^+=\{\widetilde{\alpha}|_{\mathfrak t_1^{\mathbb C}} \: |\:
\widetilde{\alpha}\in P\}.
\end{equation}
For $1 \le j \le n$, let $H_j \in \mathfrak t$ be the element with zero
diagonal entries except the $j^{th}$ entry which is $\sqrt{-1}$.
%$i=\sqrt{-1}$.
If $a \in N_G(T)$ then $aH_ja^{-1}=H_{\sigma(j)}$ for some permutation
$\sigma$ of the set $\{1,2, \ldots, n\}$, since
$Ad (a)H_j=aH_ja^{-1} \in \mathfrak t$
has the same eigenvalues of $H_j$.
One can show that the map $a \mapsto \sigma$ defines an isomorphism
$aT \mapsto \sigma$ of the Weyl group $W=N_G(T)/T$ of $U(n)$ onto the
symmetric group $S_n$ on $n$ letters such that for
$(\lambda,H)\in (\mathfrak t^{\mathbb C})^*\times
\mathfrak t^{\mathbb C}$
the action of $w \in W$,
\begin{equation}
(w\cdot \lambda)(H) \quad {\stackrel{\rm def}{=}} \quad
\lambda (Ad (a^{-1})H)=\lambda(a^{-1}Ha), \qquad w=aT,
\end{equation}
goes over to the action of $S_n$ on $(\mathfrak t^{\mathbb C})^*$
given by
\begin{equation}
(\sigma\cdot \lambda)\left( H
= {\rm diag}\,(H_1, ..., H_n)\right) =\lambda\left(
{\rm diag}\,(H_{\sigma(1)}, ..., H_{\sigma(n)})\right)
\mbox{.}
\end{equation}
%\begin{equation}  \label{6.6}
%(\sigma\cdot \lambda)\left( H
%=\begin{pmatrix} H_1 & & \\ & \ddots & \\ & & H_n \end{pmatrix}\right)
%=\lambda\left(
%\begin{pmatrix} H_{\sigma(1)} & & \\ & \ddots & \\ & & H_{\sigma(n)}
%\end{pmatrix} \right).
%\end{equation}

Sometimes it is convenient to change signs and work with
$\sigma^- {\stackrel{\rm def}{=}} -\sigma$ in place of $\sigma$.
Then $J$, $\omega_\sigma$ are replaced by
$J^-=-J$, $\omega_{\sigma^-}=(-1)^n\omega_\sigma$.
Formula (\ref{6.5}) for the choice $c=-\sqrt{-1}$ then assumes the form
\begin{equation}  \label{6.7}
\int_{G/T}e^{\sqrt{-1} \langle Ad (a^{-1})x,x_0 \rangle} \cdot
\frac{\omega_{\sigma^-}(aT)}{(2\pi)^n}
= \sum_{w\in W}\frac{e^{\sqrt{-1}(w\cdot \lambda_{x_0})(x)}}
{\prod_{\beta \in P}(w \cdot \beta)(x)}
\end{equation}
for regular elements $x, x_0 \in \mathfrak t$,
but where we now assume that
$i\beta(x_0)\stackrel{\eta}{<}0$ for all $\beta \in P$.
That is, for
\begin{equation}  \label{6.8}
x_0 = {\rm diag}\,\left(\sqrt{-1}t_1, ..., \sqrt{-1}t_n\right),
\,\,\,\,\,\,\,
x = {\rm diag}\,\left(\sqrt{-1}\theta_1, ..., \sqrt{-1}\theta_n\right)
\in \mathfrak t
\end{equation}
%\begin{equation}  \label{6.8}
%x_0 = \begin{pmatrix}it_1 & & 0\\ & \ddots & \\ 0 & &
%it_n \end{pmatrix},
%\qquad
%x =\begin{pmatrix}i\theta_1 & & 0\\ & \ddots & \\ 0 & & i\theta_n
%\end{pmatrix}
%\in \mathfrak t
%\end{equation}
the regularity condition is that the diagonal entries are all distinct,
and condition $\eta$ is that $t_1>t_2> \ldots >t_n$.
If we write
$\sigma(x) {\stackrel{\rm def}{=}}
{\rm diag}\,(\sqrt{-1}\theta_{\sigma(1)}, ...,
\sqrt{-1}\theta_{\sigma(n)})$
%$\sigma(x) {\stackrel{\rm def}{=}}
%\begin{pmatrix}
%i\theta_{\sigma(1)} & & 0\\ & \ddots & \\ 0 & & i\theta_{\sigma(n)}
%\end{pmatrix}$
for $\sigma \in S_n$, then for
\begin{equation}
P(n) \quad {\stackrel{\rm def}{=}} \quad
\frac{1}{2}n(n-1)=\frac 12 \dim U(n)/T,
\end{equation}
we obtain from the above remarks and data for $U(n)$ its
localization formula:
\begin{eqnarray}  \label{6.9}
&& \int_{U(n)/T}e^{-\sqrt{-1} \Tr (xax_0a^{-1})} \cdot
\frac{\stackrel{\rm P(n)-times}
{\overbrace{\sigma^- \wedge \ldots \wedge \sigma^-}}}
{P(n)!(2\pi)^{P(n)}}(aT)
\nonumber \\
& = & \frac{1}{(-1)^{P(n)/2}}\frac{1}{\prod_{r<s}(\theta_r-\theta_s)}
\sum_{\sigma \in S_n} (\mbox{sgn} \sigma) e^{-\sqrt{-1}
\Tr (\sigma(x)x_0) }
\nonumber \\
& = & \frac{1}{(-1)^{P(n)/2}}\frac{1}{\prod_{r<s}(\theta_r-\theta_s)}
\det \left[ e^{\sqrt{-1}\theta_i t_j}\right]
\end{eqnarray}
by (\ref{6.7}), where we have used that
\begin{equation}
\prod_{r<s}(\theta_\sigma(r)-\theta_\sigma(s))
=(\mbox{sgn}\, \sigma)\prod_{r<s}(\theta_r-\theta_s)
\end{equation}
and the expansion
$\det A = \sum_{\sigma \in S_n} (\mbox{sgn}\, \sigma) A_{1 \sigma(1)}
\ldots A_{n\sigma(n)}$
of a determinant.

Formula (\ref{6.9}) leads to the Itzykson-Zuber formula as we
now indicate.
Define a measure $\mu_0$ on $G$ by
\begin{equation}  \label{6.10}
\int_G f(a) \,d\mu_0(a) = \int_{G/T} \left[ \int_T f(at)\,dt\right]
\frac{\omega_\sigma^- (aT)}{(2\pi)^{\frac 12 \dim\, G/T}}
\end{equation}
where $dt$ denotes normalized Haar measure on $T$: $\int_T 1 \,dt=1$;
$f$ is any continuous function on $G$.
From the $G$-invariance of $\sigma$
(see the remarks following Theorem \ref{sympl})
it follows that $\omega_{\sigma^-}$ is also $G$-invariant and
that $\mu_0$
is therefore a Haar measure on $G$.
The choice $f=1$ gives $\int_G 1 \,d\mu_0=v(G/T)$ for
\begin{equation}  \label{6.11}
v(G/T) \quad {\stackrel{\rm def}{=}} \quad
\int_{G/T} \frac{\omega_{\sigma^-}}{(2\pi)^{\frac 12 \dim G/T}},
\end{equation}
which means that
$\mu {\stackrel{\rm def}{=}} {\mu _0}[v(G/T)]^{-1}$
is normalized Haar measure on $G$.
To compute $v(G/T)$ for $G=U(n)$ choose $\theta _j=\varepsilon (n-j)$
in (\ref{6.8}) for $x$.
Then the determinant in (\ref{6.9}) is Vandermonde's determinant
$=\prod_{r<s}(e^{\sqrt{-1}\varepsilon t_r}-
e^{\sqrt{-1}\varepsilon t_s})$,
and $\theta _r-\theta _s=\varepsilon (s-r)$.
The right hand side of formula (\ref{6.9}) becomes
\begin{equation}
\frac{1}{(-1)^{P(n)/2}}\frac{1}{\prod_{r<s}(s-r)}\prod_{r<s}
\frac{(e^{\sqrt{-1}\varepsilon t_r}-e^{\sqrt{-1}\varepsilon t_s})}
{\varepsilon}
\end{equation}
whose limit as $\varepsilon \rightarrow 0$ is (by L'Hospital's rule)
$\prod_{r<s}\frac{t_r-t_s}{s-r}$.
On the other hand, the limit as $\varepsilon \rightarrow 0$
of the left hand side of formula (\ref{6.9}) is $v(U(n)/T)$,
since $x\rightarrow 0$ as $\varepsilon \rightarrow 0$.
Thus we see that
(for $P(n) {\stackrel{\rm def}{=}} \frac{1}{2}n(n-1)
=\frac 12 \dim U(n)/T$)
\begin{equation}  \label{6.12}
v(U(n)/T) \quad {\stackrel{\rm def}{=}} \quad
\int_{U(n)/T} \frac{\omega_{\sigma^-}}{(2\pi)^{P(n)}}
= \prod_{r<s} \frac{t_r-t_s}{s-r}
= \frac{\prod_{r<s}(t_r-t_s)}{\prod_{k=0}^{n-1}k!}
\end{equation}
is the symplectic volume of the co-adjoint orbit defined
by $x_0$ in (\ref{6.8}).
Finally, in (\ref{6.9}) choose
$f(a) {\stackrel{\rm def}{=}} e^{-\sqrt{-1}
\Tr (xax_0 a^{-1})}$ for $x, x_0$ in
(\ref{6.8}).  Then for $b\in T$, $f(ab)=f(a)$ as $b$ commutes with $x_0$.
Keeping in mind that
$\mu {\stackrel{\rm def}{=}} {\mu_0}[v(U(n)/T)]^{-1}$
is normalized Haar measure on $U(n)$,
we obtain from equations (\ref{6.9}), (\ref{6.10}), (\ref{6.12})
\begin{eqnarray}  \label{6.13}
\!\!\!\!\!
\int_{U(n)}e^{-\sqrt{-1}\Tr (xax_0a^{-1})} \,d\mu(a)
& = & \frac{1}{v(U(n)/T)}\int_{U(n)/T}e^{-\sqrt{-1}\Tr (xax_0a^{-1})}
\cdot \frac{\omega_{\sigma^-}(aT)}{(2\pi)^{P(n)}}
\nonumber \\
\nonumber \\
& = & \frac{(\prod_{k=0}^{n-1}k!)i^{-P(n)}}{\prod_{r<s}(t_r-t_s)}
\frac{1}{\prod_{r<s}(\theta_r-\theta_s)} \det \left[
e^{e^{\sqrt{-1}}\theta_it_j}\right]\!,
\end{eqnarray}
which is the Itzykson-Zuber formula \cite{ItzyksonZuber}, obtained here
under the very general umbrella of localized equivariant cohomology --
and under the above assumptions that $x, x_0$ in (\ref{6.8}) satisfy
$t_1>t_2>\ldots>t_n$ with $\theta_r \ne \theta_s$ for $r\ne s$.

Formula (\ref{6.13}) can be formulated in the general context of an
arbitrary
compact, connected Lie group $G$ that we have been considering.
For this we need a formula for $v(G/T)$ which replaces that in (\ref{6.12})
for $U(n)$. It is given as follows.
Given $\beta \in P$ there is a unique element
$H_\beta \in \sqrt{-1} \mathfrak t_1$
such that for every $H \in t^{\mathbb C}$,
$\beta(H)= \langle H,H_\beta \rangle$, where $\langle \cdot,\cdot \rangle$
denotes the inner product on $\mathfrak g^{\mathbb C}$ that naturally
extends
$\langle \cdot,\cdot \rangle$ on $\mathfrak g$.
If $2 \delta_P {\stackrel{\rm def}{=}} \sum_{\beta \in P} \beta$, then
\begin{equation}  \label{6.14}
v(G/T) \quad {\stackrel{\rm def}{=}} \quad
\int_{G/T} \frac{\omega_{\sigma^-}}{(2\pi)^{\frac 12 \dim G/T}}
= \prod_{\beta\in P} \frac{\beta(-\sqrt{-1}x_0)}{\delta_P(H_\beta)},
\end{equation}
again for $x_0 \in \mathfrak t$ regular with $\sqrt{-1}\beta
(x_0)<0$ for all $\beta \in P$. For $G=U(n)$, for example, with
$x_0$ in (\ref{6.8}), we have seen that
$P=\{\widetilde{\alpha_{rs}}\,|\, 1 \leq r<s\leq n \}$ where
$\widetilde{\alpha_{rs}}(H)=H_r-H_s$ for $ H={\rm diag}
\left(H_1,..., H_n\right) \in \mathfrak t^{\mathbb C}. $
%$H=\begin{pmatrix}
%H_1 &        &     \\
%        & \ddots &     \\
%        &        & H_n
%\end{pmatrix} \: \in \mathfrak t^{\mathbb C}$.
Therefore $\prod_{\beta \in P} \beta(-\sqrt{-1}x_0)=\prod_{r<s}(t_r-t_s)$.
Also $\langle Z,W \rangle =\Tr (Z\overline{W}^t)$ for
$Z,W \in \mathfrak g^{\mathbb C}= \mathfrak{gl}(n,{\mathbb C})$,
$H_\beta = {\rm diag}(0,...,1,...,-1,...,0)$
%$H_\beta = \begin{pmatrix}
%0 &        &   &        &   &        &   \\
%  & \ddots &   &        &   &        &   \\
%  &        & 1 &        &   &        &   \\
%  &        &   & \ddots &   &        &   \\
%  &        &   &        &-1 &            \\
%  &        &   &        &   & \ddots &   \\
%  &        &   &        &   &        & 0
%\end{pmatrix}$
for $\beta=\widetilde{\alpha_{rs}}$, where $1,-1$ appear in the $r^{th}$
and
$s^{th}$ row, and $\delta_P(H_\beta)=s-r$.
Thus in this case formula (\ref{6.14}) reduces to formula (\ref{6.12}).
Given (\ref{6.14}) one can now repeat the argument that followed
(\ref{6.12}),
where now one takes
$f(a)=e^{\sqrt{-1} \langle Ad (a^{-1})x,x_0 \rangle}$ in (\ref{6.10})
for $x\in \mathfrak t$ also regular.
For $t\in T$ we still have $f(at)=f(a)$ since
\begin{equation}
\langle Ad ((at)^{-1})x,x_0 \rangle
= \langle Ad (t^{-1}) Ad (a^{-1})x,x_0 \rangle
= \langle Ad ((a^{-1})x, Ad (t)x_0 \rangle
\end{equation}
(by the $Ad(G)-$invariance of $\langle \cdot,\cdot \rangle$)
$= \langle Ad (a^{-1})x,x_0 \rangle$, as
\begin{equation}
T=\{ b \in G |\: Ad(b)x =y \: \forall y \in \mathfrak t\}.
\end{equation}
Following the argument exactly as given for $U(n)$ one obtains
(again for $\sqrt{-1}\beta(x_0)<0$ for all $\beta \in P$)
by (\ref{6.7}), (\ref{6.10}), (\ref{6.14})
\begin{equation}
\int_G e^{\sqrt{-1} \langle Ad (a^{-1})x,x_0 \rangle} \,d\mu(a)
= \prod_{\beta\in P}\frac{\delta_P(H_\beta)}{\beta(-\sqrt{-1}x_0)}
\sum_{w\in W}\frac{e^{\sqrt{-1}(w\cdot \lambda_{x_0})(x)}}
{\prod_{\beta \in P}(w\cdot \beta)(x)},
\end{equation}
which is a formulation of the Itzykson-Zuber formula for an
arbitrary compact, co- \\
nnected Lie group $G$.

% Matvei's section
\section{Localization formula for non-compact group actions}

In this section we describe a generalization of the B-V localization
formula
to non-compact group actions. In the classical localization formula
(\ref{For1}) it was assumed that both the manifold $M$ and the group $G$
are compact.
The compactness of $M$ ensures convergence of the integral
$\int_M [\tau]$, and the compactness of $G$ implies existence of a
$G-$invariant Riemannian metric $\langle \cdot , \cdot \rangle$ on $M$
which was used in the proof of (\ref{For1}).
When $G$ is not compact such a Riemannian metric may not exist.
Now, pick an element $X \in \mathfrak g$.
Then $M$ being compact and the space $Z(M, X, s)$ = kernel of $d_{X,s}$
on $\Lambda_XM$ being non-zero together imply that the vector field $X^*$
comes from the action of some {\em compact} group $G'$ and we could apply
the B-V localization formula (\ref{For1}) to $G'$ instead of $G$.
Thus, in order to have a truly new result where the action of $G$ does not
factor through action of some compact group we must allow non-compact
manifolds.

At a first glance it appears that the formula fails when $G$ is not
compact.
For example, let us consider $G = SL(2, \mathbb R)$ and let us take an
element $f$ in the dual of the Lie algebra
$\mathfrak g = \mathfrak {sl}(2,\mathbb R)$ defined by
\begin{equation}
f:\: \begin{pmatrix} a & b \\ c & -a \end{pmatrix} \mapsto b-c.
\end{equation}
Let ${\mathcal O} \subset \mathfrak {sl}(2,\mathbb R)^*$
denote the co-adjoint orbit of $f$.
Like all co-adjoint orbits, $\mathcal O$ possesses a canonical
symplectic structure $\sigma$ which is the top degree part of the
equivariantly closed form $s^{-1}\tau^X = s^{-1}(J(X), 0, s\sigma)$.
Although $\cal O$ is not compact, the symplectic volume
$\int_{\cal O} \sigma$ still exists as a distribution on
$\mathfrak {sl}(2,\mathbb R)$. Let
$\mathfrak {sl} (2,\mathbb R)'_{split} \subset \mathfrak {sl}(2,\mathbb R)$
be the open subset consisting of $X \in \mathfrak {sl}(2,\mathbb R)$
with distinct real eigenvalues.
In other words, $\mathfrak {sl} (2,\mathbb R)'_{split}$ consists of all
elements in $\mathfrak {sl}(2,\mathbb R)$ conjugate to
${\rm diag}(\lambda, -\lambda)$,
%$\begin{pmatrix} \lambda & 0 \\ 0 & -\lambda \end{pmatrix}$,
for some $\lambda \in \mathbb R - \{0\}$.
Now, if we take any element $X \in \mathfrak {sl} (2,\mathbb R)'_{split}$,
then one can see that the vector field $X^*$ on $\cal O$ generated by $X$
has no zeroes.
Thus, if there were a fixed point integral localization
formula like in the case of compact groups,
this formula would suggest that the distribution
determined by $\int_{\cal O} \sigma$ vanishes on the open set
$\mathfrak {sl} (2,\mathbb R)'_{split}$.
But it is known that the restriction of $\int_{\cal O} \sigma$ to
$\mathfrak {sl} (2,\mathbb R)'_{split}$ is {\em not} zero.

On the other hand, recent results from representation theory, namely
the two character formulas for representations of reductive Lie groups
due to M.~Kashiwara, W.~Rossmann, W.~Schmid and K.~Vilonen described in
\cite{Sch}, \cite{SchV2} strongly suggest that the B-V localization
formula should extend to actions of non-compact groups.
Heuristically, the failure of the localization formula in the above example
can be attributed to the lack of zeroes of $X^*$, as if they
``ran away to infinity.''

This discussion demonstrates two immediate challenges to having
a localization formula when the acting group $G$ is not compact.
First of all we must allow non-compact manifolds $M$ or homology cycles
with infinite support.
But then we need to worry about convergence of the integral $\int_M \tau$.
We will resolve this problem by restricting the class of forms that we will
integrate and by introducing a new (weaker) notion of convergence of
integrals in the sense of distributions on $\mathfrak g$.
Secondly, for arbitrary non-compact manifolds or cycles with infinite support,
the zeroes of $X^*$ tend to ``run away to infinity.''
Since we cannot have a localization formula for all manifolds and cycles,
we will specify a class of cycles for which all zeroes of $X^*$ are
accounted for and the localization formula holds.

The statement of the new localization formula uses the language of
algebraic geometry.
We consider pairs of Lie groups: a real group $G$ sitting inside
a complex one $G_{\mathbb C}$.
For example:
\begin{equation}
\begin{matrix}
GL(n, \mathbb R) & \subset & GL(n, \mathbb C) \\
GL^+(n, \mathbb R) & \subset & GL(n, \mathbb C) \\
U(n) & \subset & GL(n, \mathbb C)
\end{matrix}
\qquad
\begin{matrix}
SL(n, \mathbb R) & \subset & SL(n, \mathbb C) \\
SO(n) & \subset & SL(n, \mathbb C)  \\
SU(n) & \subset & SL(n, \mathbb C)
\end{matrix}
\qquad
Sp(n, \mathbb R) \subset Sp(n, \mathbb C)
\end{equation}
More precisely, we fix a connected complex algebraic linear reductive Lie
group $G_{\mathbb C}$ which is defined over $\mathbb R$.
We will be primarily interested in a real Lie subgroup
$G \subset G_{\mathbb C}$ lying between the group of real points
$G_{\mathbb C}(\mathbb R)$ and the identity component
$G_{\mathbb C}(\mathbb R)^0$.
%We regard $G$ as a real reductive Lie group.

Our ambient space will be the holomorphic cotangent space $T^*M$
of a smooth
complex projective variety $M$ on which $G_{\mathbb C}$ acts algebraically.
We will also assume that the maximal complex torus
$T_{\mathbb C} \subset G_{\mathbb C}$
(i.e. the maximal abelian subgroup of $G_{\mathbb C}$
isomorphic to the product of several copies of $\mathbb C - \{0\}$)
acts on $M$ with isolated fixed points.
Then there are only finitely many fixed points because $M$ is compact.
%(This condition is satisfied in all applications we have in mind.)
Let $\sigma$ denote the canonical complex algebraic holomorphic symplectic
form on $T^*M$.

For a real closed submanifold $N \subset M$, we define the real
conormal space
\begin{equation}
T^*_NM = \{ \xi \in T^*M |\: \operatorname{Re}\xi |_{T_N} = 0 \}.
\end{equation}
The homology cycles over which we will integrate equivariant forms
will include
the real conormal spaces $T^*_NM$ associated to real closed $G-$invariant
submanifolds $N \subset M$ and equipped with some orientation.
An interesting example is
$G = GL(n,\mathbb R) \subset GL(n, \mathbb C) = G_{\mathbb C}$
acting naturally on a complex Grassmanian $Gr_{\mathbb C}(k,n)$.
Let $N$ be the real Grassmanian $Gr_{\mathbb R}(k,n)$ sitting inside
$Gr_{\mathbb C}(k,n)$ and the homology cycle
$C = T^*_{Gr_{\mathbb R}(k,n)} Gr_{\mathbb C}(k,n)$.

The ordinary homology cycle is defined as a {\em finite} sum of simplices
which has no boundary. Here we will consider chains $C$ which are possibly
{\em infinite} sums of simplices. In order to be able to compute the
boundary
of $C$ we require that every point in the ambient space has an open
neighborhood which intersects only finitely many simplices.
Then the boundary $\partial C$ makes sense and we say that a chain is a
(Borel-Moore homology) cycle if $\partial C=0$.
We denote by $|C|$ the support of $C$.
The Borel-Moore cycles $C \subset T^*M$ over which we will integrate
will be subject to the following three properties:
\begin{itemize}
\item
$C$ is $G-$invariant;
\item
$C$ is real Lagrangian, i.e. $\operatorname{Re} \sigma |_C \equiv 0$
and $\dim_{\mathbb R} C = \dim_{\mathbb R} M$;
\item
$C$ is {\em conic}, i.e. invariant under the scaling action of
positive reals $\mathbb R^{>0}$ on $T^*M$ (but not necessarily under
the actions of $\mathbb C - \{0\}$ or $\mathbb R - \{0\}$).
\end{itemize}
Intuitively, these cycles $C$ consist of portions of real conormal spaces
which piece together so that there is no boundary left.

\begin{example}  \label{example}  {\em
The first non-trivial example comes from the action of
$G_{\mathbb C} = SL(2,\mathbb C)$ on the projective space
$M = \mathbb CP^1$ by projective transformations.
(Recall that $\mathbb CP^1$ is diffeomorphic to the 2-sphere;
also $n = \frac 12 \dim_{\mathbb R} \mathbb CP^1 =1$.)
The group $G = SL(2,\mathbb R)$ acts on $\mathbb CP^1$ with
exactly three different orbits:
two open hemispheres and one circle which is their common boundary.
We can position $\mathbb CP^1$ in space so that the Eastern and Western
hemispheres are stable under the $SL(2,\mathbb R)-$action.
Let $H$ denote one of these two open hemispheres, say, the Western one;
and let $S^1 \subset \mathbb CP^1$ denote the circle
containing the Greenwich meridian; $S^1 = \partial H$.

\DOUBLEFIGURE{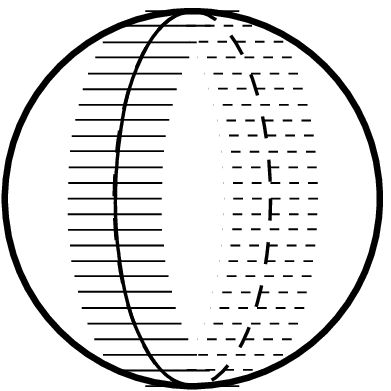} {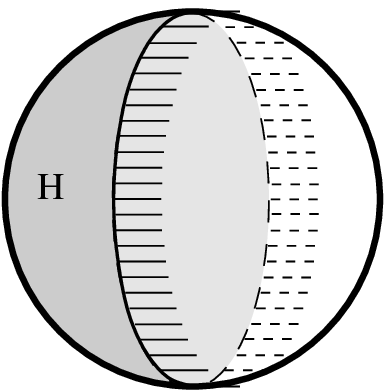}
{Cycle $T^*_{S^1} \mathbb CP^1$.}{Cycle $C_H$.}

One possible choice of the cycle $C$ is the real conormal space
$T^*_{S^1} \mathbb CP^1$ equipped with some orientation (Figure 1).
Another interesting choice of $C$ is the cycle
\begin{multline}
C_H = H \bigcup
\{ \xi \in T^*_{S^1} \mathbb CP^1 |\:
\operatorname{Re} \langle \xi, v \rangle \ge 0\,,  \\
\text{ for all tangent vectors $v$ pointing outside of $H$} \}
\end{multline}
(Figure 2).
Its orientation is determined by the orientation on $\mathbb CP^1$
which induces the orientation on $H$ and which in turn determines the
orientation on all of $C_H$.
$\Box$
}\end{example}

Let $U$ be a maximal compact subgroup of $G_{\mathbb C}$.
For instance, if $G_{\mathbb C}$ is $GL(n, \mathbb C)$ or
$SL(n, \mathbb C)$ one can take the group $U$ to be
$U(n)$ or $SU(n)$ respectively.
Then, letting $\mathfrak u$ and $\mathfrak g_{\mathbb C}$
denote the Lie algebras of $U$ and $G_{\mathbb C}$ respectively,
we have an isomorphism
$\mathfrak u \otimes_{\mathbb R} \mathbb C \simeq \mathfrak g_{\mathbb C}$.
We denote by $\Lambda^{(p,q)}M$ the space of complex-valued
differential forms of type $(p,q)$ on $M$,
that is the space of forms which are $p-$holomorphic and
$q-$antiholomorphic.
Recall that $n = \frac 12 \dim_{\mathbb R}M$.
We consider forms
$\alpha^X = (\alpha_0^X, \alpha_1^X, \dots, \alpha_{2n}^X) \in \Lambda M$
depending on $X \in \mathfrak g_{\mathbb C}$
and which satisfy the following three conditions:
\begin{itemize}
\item
The assignment $X \mapsto \alpha^X \in \Lambda M$
depends holomorphically on $X \in \mathfrak g_{\mathbb C}$;
\item
For each $k \in \mathbb N$ and each $X \in \mathfrak g_{\mathbb C}$,
$$
\alpha_{2k}^X \in
\bigoplus_{\begin{matrix} p+q=2k \\ p \ge q \end{matrix}} \Lambda^{(p,q)}M;
$$
\item
For each $X \in \mathfrak u \subset \mathfrak g_{\mathbb C}$,
we have $\alpha^X \in Z(M, X, s)$, i.e.
$$
d_{X,s} \alpha^X =0
\qquad \text{and} \qquad \theta(X)\alpha^X =0.
$$
\end{itemize}

\begin{example} {\em
A $U-$equivariant characteristic form $\alpha^X \in \Lambda M$
depending on $X \in \mathfrak u$
associated to a $U$-equivariant vector bundle over $M$
(see Section 7.1 of \cite{berline4}) satisfies the third condition.
Since it depends on $X \in \mathfrak u$ polynomially,
$\alpha^X$ extends uniquely from $\mathfrak u$ to $\mathfrak g_{\mathbb C}$
so that the first condition is satisfied.
Finally, for each $X \in \mathfrak g_{\mathbb C}$,
\begin{equation}
\alpha^X \in \bigoplus_k \Lambda^{(k,k)}M,
\end{equation}
so that the second condition is satisfied too.
This is the most important class of forms satisfying these conditions.
$\Box$
}\end{example}

Let $J(X): T^*M \to \mathbb C$ be the ordinary moment map:
\begin{equation}
J(X): \xi \mapsto \langle \xi, X^* \rangle,
\qquad \xi \in T^*M, \: X \in \mathfrak g_{\mathbb C}.
\end{equation}
The integrals will be defined as distributions on $\mathfrak g$,
so let $\varphi \in {\cal C}_c^{\infty}(\mathfrak g)$ be a test function,
and let $dX$ denote the Lebesgue measure on $\mathfrak g$.
The new localization formula will apply to integrals
of the following kind:
\begin{equation}  \label{int}
\int_C \Bigl( \int_{\mathfrak g} e^{J(X)(\xi) + \sigma}
\wedge \varphi(X) \alpha^X \,dX \Bigr)_{2n},
\qquad X \in \mathfrak g ,\: \xi \in |C| \subset T^*M.
\end{equation}
The inside integral
$
\int_{\mathfrak g} e^{J(X)(\xi) + \sigma} \wedge \varphi(X) \alpha^X \,dX
$
is essentially the Fourier transform of $\varphi(X) \alpha^X$ which decays
rapidly in the imaginary directions of
$\mathfrak g_{\mathbb C}^* \simeq \mathfrak g^* \oplus i \mathfrak g^*$.
We denote by
$
\operatorname{supp}(\sigma|_{C})
$
the closure in $T^*M$ of the set of smooth points of the support
$|C|$ where $\sigma|_{|C|} \ne 0$.
Then integral (\ref{int}) converges provided that the moment map $J$,
regarded as a map
\begin{equation}
J: \: T^*M \ni \xi \mapsto J(\cdot)(\xi) \in \mathfrak g_{\mathbb C}^*,
\end{equation}
is {\em proper} on $\operatorname{supp}(\sigma|_C)$
(meaning that the $J$-preimage of every compact set in
$\mathfrak g_{\mathbb C}^*$ is compact in
$\operatorname{supp}(\sigma|_C)$).
In particular, (\ref{int}) is well-defined when $J$ is proper on $|C|$.

Now the main result of \cite{L4} says that if the support of $\varphi$
lies in $\mathfrak g'$ ($\mathfrak g$ without a finite number of
certain hypersurfaces) then the integral (\ref{int}) can be rewritten as
\begin{equation}
\int_C \Bigl( \int_{\mathfrak g}
e^{J(X)(\xi) + \sigma} \wedge \varphi(X) \alpha^X \,dX
\Bigr)_{2n}
= \int_{\mathfrak g} F_{\alpha}(X) \varphi(X) \,dX,
\end{equation}
where $F_{\alpha}$ is a function on $\mathfrak g'$ given by the formula
\begin{equation}  \label{mainequation}
F_{\alpha}(X) = (-2\pi s)^n
\sum_{\scriptstyle p\in M,
\atop\scriptstyle p=\,{\rm a\,\, zero\,\, of}\,\,X^{*}}
m_p(X) \frac {\alpha_0^X(p)}{[{\rm det}{\frak L}_p(X^{*})]^{1/2}},
\end{equation}
and each $m_p(X)$ is a certain integer multiplicity.
The function $F_{\alpha}$ is invariant under the action of $G \cap U$
obtained by restricting the adjoint action of $G$ on $\mathfrak g$.

\begin{remark}
Perhaps the most striking new feature of this localization formula is the
presence of integer multiplicities $m_p(X)$'s.
Each multiplicity $m_p(X)$ equals the local contribution of $p$ to the
Lefschetz fixed point formula, as generalized to sheaf cohomology by
M.~Goresky and R.~MacPherson \cite{GM}.
Sheaves are a generalization of the notion of vector
bundles over a manifold. There is a recent construction due to
M.~Kashiwara which associates to each sheaf ${\cal F}$ on $M$
a cycle in $T^*M$ called the {\em characteristic cycle} of ${\cal F}$.
For example, the characteristic cycle of a vector bundle over $M$ of rank $k$
is the manifold $M$ itself regarded as a cycle in $T^*M$ and taken with
multiplicity $k$.
Any cycle $C$ satisfying the three conditions above can be realized as a
characteristic cycle $Ch({\cal F})$ of some $G$-equivariant sheaf ${\cal F}$
(\cite{KaScha}, \cite{SchV1}).
The multiplicities are determined in \cite{L4} in terms of
local cohomology of ${\cal F}$, where ${\cal F}$ is any sheaf
with characteristic cycle $Ch({\cal F}) = C$.
\end{remark}

\begin{remark}
The reason why the localization formula is stated in terms of
distributions
is that when the support of $C$ is not compact the integral
\begin{equation}
\int_C (e^{J(X)(\xi) + \sigma} \wedge \alpha^X)_{2n}
\end{equation}
practically never converges.

The set $\mathfrak g'$ is essentially the set of regular semisimple elements
of $\mathfrak g$ on which the denominators
$[{\rm det}{\frak L}_p(X^{*})]^{1/2}$ do not vanish.

In the special case when $C = M$ as oriented cycles,
$C$ is $U$-invariant, each multiplicity
$m_p(X)$ equals 1 and this theorem can be easily deduced from
the classical B-V localization formula (\ref{For1}).

Notice that the cycle $C$ is invariant with respect to the
action of the group $G$ which need not be compact, while the form
$\alpha^X \in \Lambda M$, $X \in \mathfrak g_{\mathbb C}$,
is required to be equivariant with respect to a different group $U$ only,
and $U$ may not preserve the cycle $C$.
\end{remark}

This localization formula has many interesting applications.
The most important of them is a geometric proof of the integral character
formula for representations of real reductive Lie groups \cite{L1}.
Article \cite{L2} gives a very accessible introduction to \cite{L1}
and explains the key ideas used there by way of examples and
illustrations.

\begin{example} {\em
In the setting of Example \ref{example} we consider the group
$G_{\mathbb C} = SL(2,\mathbb C)$ acting on the projective space
$M = \mathbb CP^1$, and take $G = SL(2,\mathbb R)$.
In this situation the set $\mathfrak g'$ is the set of regular semisimple
elements
\begin{equation}
\mathfrak {sl}(2,\mathbb R)^{rs} =
\{ X \in \mathfrak {sl}(2,\mathbb R) |\:
\text{$X$ has two distinct (real or complex) eigenvalues} \}.
\end{equation}
The vector field generated by each
$X \in \mathfrak {sl}(2,\mathbb R)^{rs}$
has exactly two zeroes on $\mathbb CP^1$ located diameterally opposite
to each other.
The elements of $\mathfrak {sl}(2,\mathbb R)^{rs}$ come in two flavors.
We call an element $X \in \mathfrak {sl}(2,\mathbb R)^{rs}$ {\em elliptic}
if it has purely imaginary eigenvalues or, equivalently,
if it is conjugate
to
$\begin{pmatrix} 0 & \lambda \\ -\lambda & 0 \end{pmatrix}$
for some $\lambda \in \mathbb R - \{0\}$.
We also call an element
$X \in \mathfrak {sl}(2,\mathbb R)^{rs}$ {\em split}
if it has real eigenvalues or, equivalently, if it is conjugate
to $\begin{pmatrix} \lambda & 0 \\ 0 & -\lambda \end{pmatrix}$
for some $\lambda \in \mathbb R - \{0\}$.

\DOUBLEFIGURE {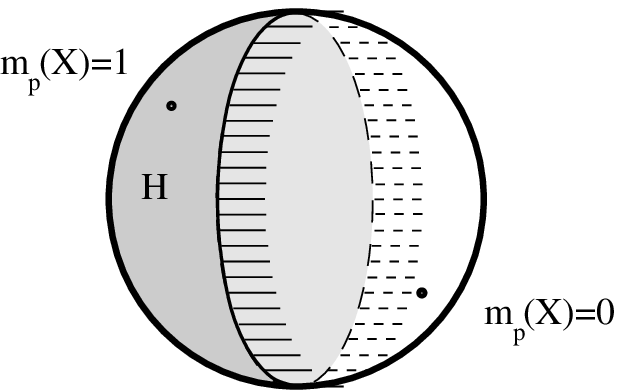}{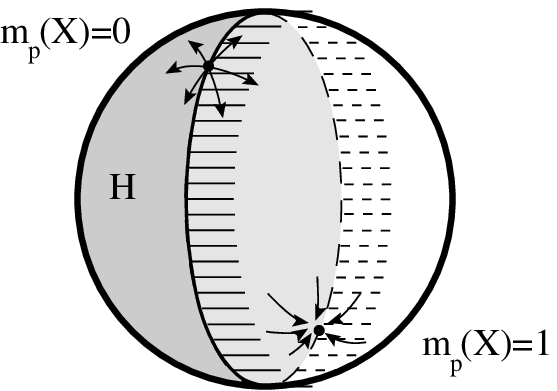}
{Zeroes of $X^*$, $X$ elliptic.}
{Zeroes of $X^*$, $X$ split.}

Consider the cycle
\begin{multline}
C_H = H \bigcup
\{ \xi \in T^*_{S^1} \mathbb CP^1 |\:
\operatorname{Re} \langle \xi, v \rangle \ge 0  \\
\text{ for all tangent vectors $v$ pointing outside of $H$} \}
\end{multline}
introduced in Example \ref{example}.
If $X \in \mathfrak {sl}(2,\mathbb R)^{rs}$ is elliptic, then the vector field
$X^*$ has one zero in the open hemisphere $H$ and the other zero in the open
hemisphere opposite to $H$.
The multiplicity
\begin{equation}
m_p(X) =
\begin{cases}
1 & \text{if $p \in H$} \\
0 & \text{if $p \notin H$}
\end{cases}
\qquad \text{$X$ is elliptic, $p$ is a zero of $X^*$}
\end{equation}
(Figure 3).
This is hardly surprising since only those zeroes of $X^*$ are expected to
make any contribution to the integral which lie in the support
of the cycle.

If $X \in \mathfrak {sl}(2,\mathbb R)^{rs}$ is split, then both zeroes of
the vector field $X^*$ lie on the boundary $S^1 = \partial H$.
While the zeroes appear to be symmetric at first, one of them counts and the
other one does not.
The symmetry is broken by the fact that one of these zeroes is stable
(the vector field $X^*$ points towards it) and the other zero is unstable
(the vector field $X^*$ points away from it).
The multiplicity
\begin{equation}
m_p(X) =
\begin{cases}
1 & \text{if $p$ is stable} \\
0 & \text{if $p$ is unstable}
\end{cases}
\qquad \text{$X$ is split, $p$ is a zero of $X^*$}
\end{equation}
(Figure 4).
This phenomenon is new and does not have analogues in compact group
actions.
$\Box$
}\end{example}

Another interesting application of the localization formula
(\ref{mainequation}) is a generalization of the Riemann-Roch-Hirzebruch
integral formula to ${\cal D}$-modules.
Its statement can be found in \cite{L5} and it uses the language of
${\cal D}$-modules (sheaves of modules over the sheaf of linear
differential operators), but its flavor can be illustrated by the
following example.

As before, $G_{\mathbb C}$ is a connected complex algebraic linear
reductive Lie group defined over $\mathbb R$ and
acting algebraically on a smooth complex projective variety $M$,
and $G \subset G_{\mathbb C}$ is a real Lie subgroup
lying between the group of real points $G_{\mathbb C}(\mathbb R)$
and the identity component $G_{\mathbb C}(\mathbb R)^0$.
Take the sheaf of sections
${\cal O}({\bf E})$ of a $G_{\mathbb C}-$equivariant algebraic line bundle
$({\bf E}, \nabla_{\bf E})$ over a $G_{\mathbb C}-$invariant open
algebraic subset
$O \subset M$ with a $G_{\mathbb C}$-invariant algebraic flat connection
$\nabla_{\bf E}$.

Let $O_{\mathbb R} \subset M$ be an open $G$-invariant subset
(which may or may not be $G_{\mathbb C}$-invariant) and
consider the cohomology spaces
\begin{equation}  \label{hspaces}
H^p(O_{\mathbb R}, {\cal O}({\bf E})).
\end{equation}
The classical Riemann-Roch-Hirzebruch formula computes the index of
${\bf E}$, i.e. the alternating sum
$
\sum_p (-1)^p \dim H^p(O_{\mathbb R}, {\cal O}({\bf E}))
$
with $O_{\mathbb R} = O = M$.
For general $O_{\mathbb R}$ and $O$, however, these dimensions can be
infinite.
To work around this problem we regard the vector spaces (\ref{hspaces})
as representations of $G$, and, as a substitute for the index, we ask
for the character of the virtual representation
$
\sum_p (-1)^p H^p(O_{\mathbb R}, {\cal O}({\bf E})).
$
(Recall that for finite-dimensional representations the value of
the character
at the identity element $e \in G$ equals the dimension of the
representation.)
This character is given by the integral formula (\ref{int}) with
concrete choices of the cycle $C$ and the form $\alpha^X$.

\end{document}